\documentclass[10pt, a4paper, fleqn]{scrartcl}
\usepackage[latin1]{inputenc}

\usepackage{a4wide}
\usepackage{amsmath}
\usepackage{amssymb}
\usepackage{amsthm}
\usepackage{paralist, calc}

\usepackage{color}              

\theoremstyle{plain}
  \newtheorem{main}{Theorem}
  \newtheorem{mcor}[main]{Corollary}
  \newtheorem{lem}{Lemma}[section]
  \newtheorem{lemma}[lem]{Lemma}
  \newtheorem{prop}[lem]{Proposition}
  \newtheorem{proposition}[lem]{Proposition}
  \newtheorem{thm}[lem]{Theorem}
  \newtheorem{cor}[lem]{Corollary}
  \newtheorem{corollary}[lem]{Corollary}

\theoremstyle{definition}
  \newtheorem{defn}[lem]{Definition}
  
  \newtheorem{rem}[lem]{Remark}
  \newtheorem{ex}[lem]{Example}

\newcommand{\F}{\ensuremath{\mathbb{F}}}    
\newcommand{\abs}[1]{\lvert#1\rvert}        
\newcommand{\mouf}{\mathbb{M}} 

\newcommand{\vct}[2]{\left(\begin{smallmatrix} {#1} \\ {#2} \end{smallmatrix}\right)} 
\newcommand{\mtr}[4]{\left(\begin{smallmatrix} {#1} & {#2} \\ {#3} & {#4} \end{smallmatrix}\right)} 

\renewcommand{\P}{\ensuremath{\mathbb{P}}}

\newcommand{\gen}[1]{\langle{#1}\rangle}

\newcommand\eps{\varepsilon}
\newcommand\sg{\sigma}

\newcommand\isomorph{\cong}

\newcommand{\SL}{{\rm SL}}
\newcommand{\PSL}{{\rm PSL}}
\newcommand{\SU}{{\rm SU}}

\newcommand{\GL}{{\rm GL}}

\newcommand{\mc}[1]{\mathcal{#1}}

\newcommand{\Gt}{G_\theta}

\DeclareMathOperator{\Aut}{Aut}
\DeclareMathOperator{\Fix}{Fix}
\DeclareMathOperator{\id}{id}

\DeclareMathOperator{\Sym}{Sym}
\DeclareMathOperator{\Str}{Str}

\numberwithin{equation}{section}

\newcommand{\Defn}[1]{\textit{#1}}

\begin{document}

\title{{Iwasawa decompositions of split Kac--Moody groups}}
\author{Tom De Medts \and Ralf Gramlich\thanks{The second author gratefully acknowledges a Heisenberg fellowship by the Deutsche Forschungsgemeinschaft.} \and Max Horn}
\maketitle



\begin{abstract}
In this article we characterize the fields over which connected split semisimple algebraic groups and split Kac-Moody groups admit an Iwasawa decomposition.
\end{abstract}

\section{Introduction}

The Iwasawa decomposition of a connected semisimple complex Lie group or a connected semisimple split real Lie group is one of the most fundamental observations of classical Lie theory. It implies that the geometry of a connected semisimple complex resp.\ split real Lie group $G$ is controlled by any maximal compact subgroup $K$. Examples are Weyl's unitarian trick in the representation theory of Lie groups, or the transitive action of $K$ on the Tits building $G/B$. In the case of the connected semisimple split real Lie group of type $G_2$ the latter implies the existence of an interesting epimorphism from the real building of type $G_2$, the real split Cayley hexagon, onto the real building of type $A_2$, the real projective plane, by means of the epimorphism $\mathrm{SO}_4(\mathbb{R}) \to \mathrm{SO}_3(\mathbb{R})$ (see \cite{Gramlich:1998}). This epimorphism cannot be described using the group of type $G_2$, because it is quasisimple.

\medskip
To be able to transfer these ideas to a broader class of groups that includes the class of split Kac--Moody groups over arbitrary fields, we extend the notion of an Iwasawa decomposition in the following way, which is inspired by the fact that a maximal compact subgroup of a semisimple split Lie group is centralized by an involution:

\begin{defn} \label{Iwasawa}
A group $G$ with a twin $BN$-pair (cf.\ Definition \ref{twin}) admits an \Defn{Iwasawa decomposition}, if there exist an involution $\theta\in\Aut(G)$ such that
\begin{enumerate}
\item $B_+^\theta = B_-$ and 
\item $G = \Gt B_+$ where $\Gt:=\Fix_G(\theta)$.
\end{enumerate}
\end{defn}

Our interest in Iwasawa decompositions stems from geometric group theory: the group $\Gt$ acts with a fundamental domain on the flag complex of the building $G/B_+$, which is simply connected by \cite[Theorem 4.127]{Abramenko/Brown:2008}, \cite[Theorem 13.32]{Tits:1974}, if its rank is sufficiently large. Hence Tits' Lemma \cite[Lemma 5]{Pasini:1985}, \cite[Corollary 1]{Tits:1986} yields a presentation of $\Gt$ by generators and relations. 

Since $\Gt$ need not be finitely generated, this presentation is usually formulated as a universal enveloping result of an amalgam.
The following theorem specifies the presentation we have in mind. Refer to Section \ref{bnrgd} for a definition of a root group datum and the construction of the twin $BN$-pair resulting from it.

\begin{main} \label{thm2}
Let $G$ be a group with a root group datum $\{U_\alpha\}_{\alpha\in\Phi}$,
and assume that $G = \Gt B_+$ is an Iwasawa decomposition of $G$ with
respect to an involution $\theta$.
Furthermore, let $\Pi$ be a system of fundamental roots of $\Phi$ and, for
$\{\alpha, \beta\} \subseteq \Pi$, let $X_{\alpha,\beta} := \gen{U_\alpha,U_{-\alpha},U_\beta,U_{-\beta}}$.

Then $\theta$ induces an involution on each $X_{\alpha,\beta}$, and $\Gt$ is the universal enveloping group of the amalgam $((X_{\alpha,\beta})_\theta)_{\{\alpha, \beta \} \subseteq \Pi}$ of the subgroups of the $\theta$-fixed elements of the groups $X_{\alpha,\beta}$.
\end{main}

The proof of Theorem \ref{thm2} can be adapted as follows from what has been done in \cite{Gloeckner/Gramlich/Hartnick} and \cite{Gramlich:2006} for compact real forms of complex Lie groups and of complex Kac-Moody groups.
The involution $\theta$ induces an involution of each group $X_{\alpha,\beta}$ by Corollary \ref{3.5}.
 By the Iwasawa decomposition the group $\Gt$ acts with a fundamental domain on the flag complex $\Delta$ associated to the building $G/B_+$. Choose $F$ to be a fundamental domain of $\Delta$ stabilized by the standard torus $T:=B_+ \cap N$ of $G$ arising from the positive $BN$-pair of $G$. The stabilizers of the simplices of $F$ of dimension zero and one with respect to the natural action of $G$ on $\Delta$ are exactly the groups $(X_\alpha)_\theta T$ and $(X_{\alpha\beta})_\theta T$, $\alpha, \beta \in \Pi$. By the simple connectedness of building geometries of rank at least three (cf.\ \cite[Theorem 4.127]{Abramenko/Brown:2008} or \cite[Theorem 13.32]{Tits:1974}) plus Tits' Lemma (cf.\ \cite[Lemma 5]{Pasini:1985}, \cite[Corollary 1]{Tits:1986}) the group $\Gt$ equals the universal enveloping group of the amalgam $((X_{\alpha\beta})_\theta T)_{\alpha, \beta \in \Pi}$. Finally, the torus $T$ equals the universal enveloping group of the amalgam $(T_{\alpha\beta})_{\alpha,\beta \in \Pi}$, where $T_{\alpha\beta}$ denotes the maximal torus $T \cap X_{\alpha\beta}$ of $X_{\alpha\beta}$, and so by \cite[Lemma 29.3]{Gorenstein/Lyons/Solomon:1995} the group $G$ actually equals the universal enveloping group of the amalgam $((X_{\alpha\beta})_\theta)_{\alpha, \beta \in \Pi}$.

\begin{rem}
For a group $G$ with a 3-spherical root group datum and an involution $\theta$ of $G$ that interchanges $B_+$ and $B_-$, the local-to-global criterion proposed in \cite{Devillers/Muehlherr:2007} allows one to study amalgam presentations of the subgroup $G_\theta$ of $\theta$-fixed elements even if $G$ does not admit an Iwasawa decomposition with respect to $\theta$. The main obstruction encountered when using that approach is establishing that certain two-dimensional simplicial complexes are simply connected. We refer the interested reader to \cite{Devillers/Gramlich/Muehlherr:2009} and \cite{Gramlich/Witzel} for an application of the criterion from \cite{Devillers/Muehlherr:2007}.
\end{rem}

\medskip
From Theorem \ref{thm2} rises the question which groups actually admit an Iwasawa decomposition. In the literature one can find a lot of information on Iwasawa decompositions for prescribed ground fields, usually over the complex numbers $\mathbb{C}$, the real numbers $\mathbb{R}$, or real closed fields, cf.\ \cite{Balan/Dorfmeister:2001}, \cite{Beltita}, \cite{Grosshans:1972}, \cite{Helgason:1978}, \cite{Kroetz}, \cite{Pati/Parashar:1998}. In \cite{Helminck/Wang:1993} on the other hand, no fixed ground field is chosen. Instead, it is shown that, under the assumption that the ground field $\mathbb{F}$ is infinite and that each square is a fourth power, the Iwasawa decomposition, the polar decomposition and the $KAK$-decomposition of the subgroup of $\mathbb{F}$-rational points of a connected reductive algebraic $\mathbb{F}$-group are equivalent. 

Neither do we choose a fixed field in this paper. Instead, in Theorem \ref{thm1} below we characterize the fields $\F$ over which a group with an $\F$-locally split root group datum (Definition \ref{rgd}) admits an Iwasawa decomposition. We point out that this class of groups contains the class of groups of $\F$-rational points of a connected split semisimple algebraic group defined over $\F$ (cf.\ \cite{Springer:1998}) and the class of split Kac--Moody groups over $\F$ (cf.\ \cite{Remy:2002}, \cite{Tits:1987}). Not every group that admits an Iwasawa decomposition admits a polar or a $KAK$-decomposition.

\medskip 
For the next definition recall that any Cartan--Chevalley involution of $(\mathrm{P})\SL_2(\F)$ is given, resp.\ induced by the transpose-inverse automorphism with respect to the choice of a basis of the natural $\SL_2(\F)$-module $\F^2$.
\begin{defn}
Let $\F$ be a field, let $\sigma$ be an automorphism of $\F$ of order $1$ or $2$, and let $G$ be a group with an $\F$-locally split root group datum $\{U_\alpha\}_{\alpha\in\Phi}$. We call an automorphism $\theta\in\Aut(G)$ a \Defn{$\sigma$-twisted Chevalley involution} of $G$, if it satisfies the following for all $\alpha\in\Phi$:
\begin{enumerate}
\item $\theta^2=\id$,
\item $U_\alpha^\theta=U_{-\alpha}$, and
\item $\theta\circ\sigma$ induces a Cartan--Chevalley involution (resp.\ its image under the canonical projection) on $X_\alpha:=\gen{U_\alpha, U_{-\alpha}} \isomorph (\mathrm{P})\SL_2(\F)$.
\end{enumerate}
\end{defn}
A $\sigma$-twisted Chevalley involution of a split Kac-Moody group can be constructed by taking the product of a sign automorphism and the field automorphism $\sigma$, see  \cite[Section 8.2]{Caprace/Muehlherr:2005}. The same is true for all split semisimple algebraic groups. Moreover, by Lemma \ref{chev-inv-existence} below a group with a $2$-spherical $\F$-locally split root group datum with $|\F| \geq 4$ also admits a $\sigma$-twisted Chevalley involution. Hence the following theorem applies to these three classes of groups.

\begin{main} \label{thm1}
Let $\F$ be a field and let $G$ be a group with an $\F$-locally split root group datum.
The group $G$ admits an Iwasawa decomposition if and only if $\F$ admits an automorphism $\sigma$ of order $1$ or $2$ such that 
\begin{enumerate}
\item $-1$ is not a norm,
\item
  \begin{enumerate}
    \item if there exists a rank one subgroup $\gen{U_\alpha,U_{-\alpha}}$
          of $G$ isomorphic to $\SL_2(\F)$, then a sum of norms is a norm, or 
    \item if each rank one subgroup $\gen{U_\alpha,U_{-\alpha}}$
          of $G$ is isomorphic to $\PSL_2(\F)$, then a sum of norms is $\pm 1$ times a norm,
  \end{enumerate}
\end{enumerate}
 with respect to the norm map $N : \F \to \Fix_\F(\sigma) : x \mapsto x x^\sigma$, and
\begin{enumerate}
\setcounter{enumi}{2}
\item $G$ admits a $\sigma$-twisted Chevalley involution.
\end{enumerate}
\end{main}
\begin{rem}
We emphasize that the above theorem does not state that each Iwasawa decomposition is realized by a $\sigma$-twisted Chevalley involution. For example, any involution $\sigma$ of $\mathbb{C}$ yields an Iwasawa decomposition $\SL_n(\mathbb{C}) = \SU_n(\mathbb{C}, \theta_\sigma) B_+$ with respect to the $\sigma$-twisted Chevalley involution $\theta_\sigma$ by the preceding theorem, because $\mathrm{Fix}_\mathbb{C}(\sigma)$ is real closed. Hence, given two distinct involutions $\sigma_1, \sigma_2 \in \mathbb{C}$, the involution $(\theta_{\sigma_1},\theta_{\sigma_2})$ yields an Iwasawa decomposition of the product $\SL_n(\mathbb{C}) \times \SL_n(\mathbb{C})$ although it is not a $\sigma$-twisted Chevalley involution.
\end{rem}

In a finite field $\mathbb{F}_{q}$ of order $q \equiv 3 \mod 4$ the element $-1$ is not a square, so that a sum of squares is plus or minus a square, because there exist exactly two square classes, cf.\ \cite[Section VI \S62]{OMeara:1973}. Therefore, by Theorem \ref{thm1}, a group with an $\F_q$-locally split root group datum where all rank one subgroups are isomorphic to $\PSL_2(\F_q)$ admits an Iwasawa decomposition with respect to the standard Chevalley involution. 
We point out that it follows from inspection of the Chevalley groups of rank two (see \cite{Steinberg:1968}) that the diagram of a group with an $\F$-locally split root group datum must be right angled (i.e., any edge of the diagram is labelled with infinity), if all its rank one groups are isomorphic to $\PSL_2(\F)$. Nontrivial such examples can be obtained using free constructions (see \cite[Example 2.8]{Caprace/Remy:2008-LN}, \cite{Remy/Ronan:2006}).

\medskip
As mentioned before, among the most prominent groups that admit $\sigma$-twisted Chevalley involutions are the connected split semisimple algebraic groups and the split Kac-Moody groups. We explicitly re-state our characterization given in Theorem \ref{thm1} for these two classes of groups.

\begin{mcor} \label{3}
Let $\F$ be a field and let $G$ be a connected split semisimple algebraic group defined over $\F$ or a split Kac-Moody group over $\F$.
The group $G(\F)$ admits an Iwasawa decomposition $G(\mathbb{F}) = \Gt(\mathbb{F}) B(\mathbb{F})$ if and only if $\F$ admits an automorphism $\sigma$ of order $1$ or $2$ such that 
\begin{enumerate}
\item $-1$ is not a norm,
\item
  \begin{enumerate}
    \item if there exists a rank one subgroup $\gen{U_\alpha,U_{-\alpha}}$
          of $G$ isomorphic to $\SL_2(\F)$, then a sum of norms is a norm, or 
    \item if each rank one subgroup $\gen{U_\alpha,U_{-\alpha}}$
          of $G$ is isomorphic to $\PSL_2(\F)$, then a sum of norms is $\pm 1$ times a norm,
  \end{enumerate}
\end{enumerate}
 with respect to the norm map $N : \F \to \Fix_\F(\sigma) : x \mapsto x x^\sigma$.
\end{mcor}

We remark that Corollary \ref{3} also holds for split reductive algebraic groups, which in the setting we chose are excluded by axiom (RGD0) of Definition \ref{rgd}. By using a slightly more general notion of a group with a root group datum (such as in \cite{Tits:1992}) one can generalize Corollary \ref{3} to the reductive case.

\subsubsection*{An application to graph theory}

Another consequence of Theorem \ref{thm1} is a combinatorial local characterization of certain graphs similar to the main result of \cite{Altmann:2007}. Before we can state that result we have to introduce some additional terminology.  

\begin{defn}
Let $\F$ be a field admitting an automorphism $\sigma$ of order $1$ or $2$ such that 
\begin{enumerate}
\item $-1$ is not a norm and 
\item a sum of norms is a norm
\end{enumerate}
with respect to the norm map $N : \F \to \Fix_\F(\sigma) : x \mapsto x x^\sigma$. Then the pair $(\F,\sg)$ is called an \Defn{Iwasawa pair}.
\end{defn}

Let $\F$ be a field with an automorphism $\sigma$ of order $2$ such that both $(\F,\sg)$ and $(\Fix_\F(\sg),\id)$ are Iwasawa pairs.
Furthermore, let $V$ be a six-dimensional $\F$-vector space endowed with an anisotropic $\sigma$-hermitian sesquilinear form. Define $\mathbf{S}(V)$ to be the graph with the $2$-dimensional subspaces of $V$ as vertices and adjacency given by orthogonality. Because of this definition it makes sense to use the symbol $\perp$ for adjacency.

There exists a well known group-theoretical characterization of $\mathbf{S}(V)$ as follows. Let $G \cong \SL_6(\mathbb{F})$, so that $G$ has an $\F$-locally split root group datum $\{ U_\alpha \}_{\alpha \in \Phi}$ of type $A_5$. Let $\theta$ be the $\sigma$-twisted Chevalley involution of the group $G$, let $K:= \Fix_G(\theta)$ and, for $\alpha \in \Phi$, denote by $(X_{\alpha})_\theta$ the fixed point subgroup of the rank one group $X_{\alpha} = \gen{U_\alpha,U_{-\alpha}}$. Then $\mathbf{S}(V)$ is isomorphic to the graph on the $K$-conjugates of $(X_{\alpha})_\theta$ with the commutation relation as adjacency. 
This setup allows us to define the graph $\mathbf{S}$ for any connected split semisimple algebraic $\F$-group, so that for instance it makes sense to use the symbol $\mathbf{S}(E_6(\F))$. In this notation, we have $\mathbf{S}(V) \cong \mathbf{S}(A_5(\F))$.

The proof of \cite[Theorem 4.1.2]{Altmann:2007}, which deals with the special case of $\mathbb{F} = \mathbb{C}$ and $\sigma$ complex conjugation,  applies verbatim to the setting just introduced, so that we obtain the following result. Recall that, if $\Delta$ is a graph, one says that a graph $\Gamma$ is \Defn{locally $\Delta$}, if for each vertex $x \in \Gamma$ the subgraph of $\Gamma$ induced on the neighbors of $x$ is isomorphic to $\Delta$.

\begin{main} \label{thm4}
Let $\F$ be a field admitting an automorphism $\sigma$ of order $2$ such that both $(\F,\sg)$ and $(\Fix_\F(\sigma),\id)$ are Iwasawa pairs. Assume that $N(\F)$ is a subset of the squares of $\Fix_\F(\sigma)$, where $N : \F \to \Fix_\F(\sigma) : x \mapsto x x^\sigma$.
Let $\Gamma$ be a connected locally $\mathbf{S}(A_5(\F))$ graph satisfying that, for every chain $x \perp w \perp y \perp z \perp x \perp y$ consisting of four distinct vertices of $\Gamma$, the vertices $w$, $x$, $z$ have $y$ as unique common neighbor if and only if the vertices $w$, $y$, $z$ have $x$ as unique common neighbor.

Then $\Gamma$ is a quotient of $\mathbf{S}(A_7(\F))$ or of $\mathbf{S}(E_6(\F))$.
\end{main}

\subsubsection*{An application to group theory}

Finally, as explained in \cite[Section~4.6]{Altmann:2007}, Theorem \ref{thm4} implies the following group-theoretic statement via a standard argument. All twisted groups are understood to be anisotropic forms of connected split semisimple algebraic groups with respect to the Chevalley involution composed with the indicated field automorphism.

\begin{main} \label{twisted}
Let $\F$ be a field admitting an automorphism $\sigma$ of order $2$ such that both $(\F,\sg)$ and $(\Fix_\F(\sigma),\id)$ are Iwasawa pairs. Assume that $N(\F)$ is a subset of the squares of $\Fix_\F(\sigma)$, where $N : \F \to \Fix_\F(\sigma) : x \mapsto x x^\sigma$.
Moreover, let $G$ be a group containing an involution $x$ and a subgroup $K \unlhd C_G(x)$ such that 
\begin{enumerate}
\item $K \cong \mathrm{SU}_{6}(\mathbb{F},\theta_\sigma)$;
\item $C_G(K)$ contains a subgroup $X \cong \SU_2(\mathbb{F},\theta_\sigma)$ with $x = Z(X)$;
\item there exists an involution $g \in G$ such that $Y := gXg$ is contained in $K$; 
\item if $V$ is a natural module for $K$, then the commutator $[Y,V] := \{ yv - v \in V \mid y \in Y, v \in V \}$ has $\mathbb{F}$-dimension $2$;
\item $G = \gen{K,gKg}$; moreover, there exists $z \in K \cap gKg$ which is a $gKg$-conjugate of $x$ and a $K$-conjugate of $gxg$.
\end{enumerate}
Then
$$G/Z(G)
\cong {\rm PSU}_{8}(\mathbb{F},\theta_\sigma) \mbox{ or } G/Z(G) \cong {}^2E_6(\mathbb{F},\theta_\sigma).$$
\end{main}

\subsubsection*{Organization of the article}

In Section \ref{bnrgd} we quickly recall the definitions of a (twin) $BN$-pair and a root group datum. In Section \ref{flipbn} we collect basic facts about flips of groups with a root group datum and provide a local-to-global argument that reduces the proof of Theorem \ref{thm1} to the study of flips of rank one groups. In Section \ref{transinv} we thoroughly study flips of $(\mathrm{P})\SL_2(\F)$ and in Section \ref{sect:moufang} we use Moufang sets in order to generalize some of our results for $(\mathrm{P})\SL_2(\F)$ to arbitrary rank one groups. In particular, we study flips of $(\mathrm{P})\SL_2(\mathbb{D})$, where $\mathbb{D}$ is an arbitrary division ring.

\subsubsection*{Acknowledgements}

The authors thank Bernhard M\"uhlherr for several discussions on the topic of this article. Moreover, the authors express their gratitude to Andreas Mars for a careful proof-reading and helpful comments on an earlier version of this paper.
They also thank the referee for further remarks, comments and suggestions.

\section{\boldmath Twin $BN$-pairs and root group data} \label{bnrgd}

\begin{defn}
We call the tuple $(G,B,N,S)$  consisting of a group $G$ with subgroups $B$ and $N$ and of a subset $S$ of the coset space $N/(B \cap N)$, a \Defn{Tits system} or a \Defn{$BN$-pair} if the following conditions are satisfied:
\begin{enumerate}
\item $G=\gen{B,N}$;
\item $T:=B\cap N$ is normal in $N$;
\item the elements of $S$ have order $2$ and generate the group $W:=N/T$, called \Defn{Weyl group};
\item $BwBsB \subset BwsB\cup BwB$ for all $w \in W$, $s \in S$;
\item $sBs \not\subset B$ for all $s \in S$.
\end{enumerate}
\end{defn}

\begin{defn} \label{twin}
Let $(G,B_+,N,S)$ and $(G,B_-,N,S)$ be Tits systems such that $B_+ \cap N=B_- \cap N$, i.e., with equal Weyl group $W$. Then $(G,B_+,B_-,N,S)$ is called a \Defn{twin Tits system} or \Defn{twin $BN$-pair} if the following conditions are satisfied, cf.\ \cite{Tits:1992}:
\begin{enumerate}
\item $B_\eps w B_{-\eps} s B_{-\eps} = B_\eps w s B_{-\eps}$ for $\eps = \pm$ and all $w \in W$, $s \in S$ such that $l(ws) < l(w)$;
\item $B_+s \cap B_- = \emptyset$ for all $s \in S$.
\end{enumerate}
The conjugates of $B_+$ and $B_-$ in $G$ are called the {\em Borel subgroups} of $G$.
\end{defn}

The twin $BN$-pair is called {\em saturated}, if $B_+ \cap B_- = B_+ \cap N = B_- \cap N$.  From a geometric point of view there is no loss in restricting one's attention to saturated twin $BN$-pairs by \cite[Lemma~6.85]{Abramenko/Brown:2008}.
A group $G$ admitting a $BN$-pair satisfies $G = \bigsqcup_{w \in W} BwB$, the \Defn{Bruhat decomposition} of $G$. For each $s \in S$ the set $P_s := B \cup BsB$ is a subgroup of $G$. A Tits system $(G, B, N, S)$ leads to a \Defn{building} whose set of chambers equals $G/B$ and whose distance function $\delta : G/B \times G/B \to W$ is given by $\delta(gB,hB) = w$ if and only if $Bh^{-1}gB = BwB$. 

A group $G$ with a twin $BN$-pair hence yields two buildings $G/B_+$ and $G/B_-$ with distance functions $\delta_+$ and $\delta_-$. Furthermore, it admits the \Defn{Birkhoff decomposition} $G = \bigsqcup_{w \in W} B_+wB_-$ from which we can define the codistance function $\delta_* : (G/B_- \times G/B_+) \cup (G/B_+ \times G/B_-) \to W$ via $\delta_*(gB_-,hB_+) = w$ if and only if $B_+h^{-1}gB_- = B_+wB_-$ and $\delta_*(hB_+,gB_-) := (\delta_*(gB_-,hB_+))^{-1}$. The tuple $((G/B_+,\delta_+),(G/B_-,\delta_-),\delta_*)$ is called the \Defn{twin building} of $G$.

In the present paper we are only interested in (twin) buildings coming from a group with a (twin) $BN$-pair. For detailed treatments of the theory of buildings we refer to \cite{Abramenko/Brown:2008}, \cite{Brown:1989}, \cite{Ronan:1989}, \cite{Tits:1974}, \cite{Weiss:2003}. Twin buildings are treated in \cite{Remy:2002}, \cite{Tits:1992}.

\medskip
Let $(W,S)$ be a Coxeter system and let $\Phi$ be the set of its roots. Following \cite{Caprace/Muehlherr:2005} a \Defn{root} is a set of the form $w \alpha_s$, where $w \in W$ and $\alpha_s = \{ x \in W \mid l(sx) = l(x) + 1 \}$ with length function $l$ of $W$ with respect to the generating set $S$. Moreover, let $\Pi$ be a system of fundamental roots of $\Phi$ and, for $\eps = \pm$, let $\Phi_\eps$ denote the set of positive, resp.\ negative roots of $\Phi$ with respect to $\Pi$. For a root $\alpha\in\Phi$, denote by $s_\alpha$ the reflection of $W$ which permutes $\alpha$ and $-\alpha$. 
For each $w \in W$, define $\Phi_w := \{ \alpha \in \Phi_+ \mid w \alpha \in \Phi_- \}$.
A pair $\{\alpha,\beta\}$ of roots is called \Defn{prenilpotent} if $\alpha\cap \beta$ and $(-\alpha)\cap (-\beta)$ are both non-empty. In that case denote by $[\alpha,\beta]$ the set of all roots $\gamma$ of $\Phi$ such that $\alpha\cap \beta\subseteq\gamma$ and $(-\alpha) \cap (-\beta) \subseteq-\gamma$, and set $]\alpha,\beta[\;:=[\alpha,\beta]\setminus\{\alpha,\beta\}$.

\medskip
The following definition of a root group datum is taken from \cite{Caprace/Remy:2008-LN}.

\begin{defn} \label{rgd}
A \Defn{root group datum} of type $(W,S)$ for a group $G$ is a family $\{U_\alpha\}_{\alpha\in\Phi}$ of subgroups (the \Defn{root subgroups}) of $G$ satisfying the following axioms, where
$U_+:=\gen{U_\alpha\mid \alpha\in\Phi_+}$ and $U_-:=\gen{U_\alpha\mid \alpha\in\Phi_-}$:
\begin{description}
\item[(RGD0)] For each $\alpha\in\Phi$, we have $U_\alpha\neq1$; moreover, $G = \gen{U_\alpha \mid \alpha \in \Phi}$.
\item[(RGD1)] For each $\alpha\in\Pi$ we have $U_\alpha \not\subseteq U_-$.
\item[(RGD2)] For each $\alpha\in \Pi$ and $u\in U_\alpha\setminus\{1\}$, there exist elements $u',u''$ of $U_{-\alpha}$ such that the product $\mu(u):=u'uu''$ conjugates $U_\beta$ onto $U_{s_\alpha(\beta)}$ for each $\beta\in\Phi$.
\item[(RGD3)] For each prenilpotent pair $\{\alpha,\beta\}\subset\Phi$, we have $[U_\alpha,U_\beta]\subset\gen{U_\gamma\mid \gamma \in ]\alpha,\beta[}$.
\item[(RGD4)] For each $\alpha \in \Pi$ there exists $\alpha' \in \Phi_{s_\alpha}$ such that $U_\beta \subseteq U_{\alpha'}$ for each $\beta \in \Phi_{s_\alpha}$. 
\end{description}

Defining 
\begin{eqnarray*}
T & := & \gen{\mu(u)\mu(v) \mid u, v \in U_\alpha \backslash \{ 1 \}, \alpha \in \Pi }, \\
N & := & \gen{\mu(u) \mid u \in U_\alpha \backslash \{ 1 \}, \alpha \in \Pi }, \\
B_+ & := & T.U_+, \\
B_- & := & T.U_-, 
\end{eqnarray*}
we obtain a {twin $BN$-pair} of $G$, which by the above leads to a twin building on which $G$ acts.
We set $X_\alpha:=\gen{U_\alpha, U_{-\alpha}}$ and $X_{\alpha,\beta}:=\gen{X_\alpha,X_\beta}$. A root group datum is called \Defn{locally split} if the group $T$ is abelian and if for each $\alpha \in \Phi$ there is a field $\F_\alpha$ such that the root group datum $\{ U_\alpha, U_{-\alpha} \}$ of $X_\alpha$ of type $A_1$ is isomorphic to the natural root group datum of $\SL_2(\F_\alpha)$ or $\PSL_2(\F_\alpha)$. A locally split root group datum is called \Defn{$\F$-locally split} if $\F_\alpha \cong \F$ for all $\alpha \in \Phi$.
\end{defn}

\section{\boldmath From local to global and back} \label{flipbn} \label{localcrit3}

We start this section by discussing very basic properties of flips of twin buildings and of groups with a twin $BN$-pair. For a more thorough treatment beyond what is needed in the present article we refer the reader to \cite{Gramlich/Horn/Muehlherr}.

\begin{defn} \label{BN-flip}
Let $G$ be a group with a saturated twin $BN$-pair $B_+$, $B_-$, $N$.
An automorphism $\theta$ of $G$ is called a \Defn{$BN$-flip}, if the following holds:
\begin{enumerate}
\item $\theta^2 = \id$;
\item $B_+^\theta = B_-$;
\item $\theta$ induces a trivial action on the Weyl group $N/T$.
\end{enumerate}
\end{defn}

\begin{rem}
We stress that the key issue in item (iii) of the preceding definition is not that $\theta$ acts on the Weyl group, but that it centralizes it; cf.\ \cite[Proposition~3.1, Remark~3.3]{Gramlich/Horn/Muehlherr}.  
\end{rem}

We now give the definition of a flip of a twin building and describe the correspondence between a $BN$-flip and a flip of the twin building induced by the twin $BN$-pair. We refer to \cite{Gramlich/Horn/Muehlherr} and \cite{Horn:2008} for a more general treatment of flips.
Originally, the concept of a flip of a twin building has been introduced in \cite{Bennett/Gramlich/Hoffman/Shpectorov:2003}.

\begin{defn}
Let $\mc{B}=(\mc{B}_+,
\mc{B}_-,\delta_{\ast})= ((\mc{C}_+,\delta_+),(\mc{C}_-,\delta_-),\delta_*) $ be a twin building.
A \Defn{building flip} is an involutory permutation $\theta$ of $\mc{C}_+ \cup \mc{C}_-$ with the following properties:
\begin{enumerate}
\item $\mc{C}_+^\theta = \mc{C}_-$;
\item $\theta$ flips the distances, i.e., for $\eps = \pm$ and for all $x, y \in \mathcal{C}_\eps$ we have $\delta_\eps(x,y)= \delta_{-\eps}(x^\theta, y^\theta)$; and
\item $\theta$ preserves the codistance, i.e., for $\eps = \pm$ and for all $x \in \mathcal{C}_\eps$, $y \in \mathcal{C}_{-\eps}$ we have $\delta_\ast (x,y) =
\delta_\ast (x^\theta, y^\theta)$.
\end{enumerate}
If, additionally,
\begin{enumerate}
\setcounter{enumi}{3}
\item there exists a chamber $c
\in \mc{C}_{\pm}$
such that $\delta_\ast(c,c^\theta)=1_W$,
\end{enumerate}
the building flip $\theta$ is called a
\Defn{Phan involution}.
\end{defn}

The following proposition, which is is a special case of \cite[Proposition 3.1]{Gramlich/Horn/Muehlherr}, shows that a $BN$-flip induces a building flip (even a Phan involution) justifying the choice of name. 
Conversely, as Bernhard Mühlherr pointed out to us, a building flip induces a $BN$-flip, if the group $G$ admits a root group datum and is center free (see \cite[Theorem~2.2.2]{Horn:2008}).

\begin{prop} \label{equivalence}
Let $G$ be a group with a saturated twin $BN$-pair inducing the twin building $\mc{B}$. Then any $BN$-flip $\theta$ of $G$ induces a Phan involution of $\mc{B}$.
\end{prop}

\begin{proof}
Recall from Definition \ref{twin} that $\mc{B}$ consists of the buildings $G/B_\eps$ with distance functions $\delta_\eps : G/B_\eps \times G/B_\eps \to W$ satisfying $\delta_\eps (gB_\eps,hB_\eps) = w$ if and only if $B_\eps g^{-1}hB_\eps = B_\eps w B_\eps$ for $\eps = \pm$. These buildings are twinned by the codistance function $\delta_* : (G/B_+ \times G/B_-) \cup (G/B_- \times G/B_+) \to W$ satisfying $\delta_\ast (gB_\eps,hB_{-\eps}) = w$ if and only if $B_\eps g^{-1}hB_{-\eps} = B_\eps w B_{- \eps}$.
By definition $\theta^2 = \id$. Moreover, $B_+^\theta = B_-$ implies that $\theta$ interchanges the two parts of the twin building. The image of $g^{-1}h \in B_\eps w B_\eps$ under $\theta$ satisfies $(g^{-1})^\theta h^\theta \in B^\theta_{\eps} w^\theta B^\theta_{\eps} = B_{- \eps} w B_{- \eps}$. Therefore $\delta_{\eps} (g B_{\eps},h B_{\eps}) = w$ implies $\delta_{- \eps} (g^\theta B_{- \eps},h^\theta B_{- \eps}) = w$, whence $\theta$ flips the distances. For $g^{-1} h \in B_\eps w B_{- \eps}$ we have $(g^{-1})^\theta h^\theta \in B^\theta_{\eps} w^\theta B^\theta_{-\eps} = B_{- \eps} w B_{\eps}$, so $\delta_{\ast} (g B_{\eps},h B_{-\eps}) = w$ implies $\delta_{\ast} (g^\theta B_{- \eps},h^\theta B_{\eps}) = w$, whence $\theta$ preserves the codistance. Finally, the chamber $B_+$ is mapped onto its opposite chamber $B_-$. Altogether, this implies that $\theta$ induces a Phan involution on the twin building $\mc{B}$.
\end{proof}


\begin{proposition} \label{prop:rgd-flip}
Let $G$ be a group with a saturated twin $BN$-pair $B_+$, $B_-$, $N$, and let $\theta$ be an automorphism of $G$ satisfying
\begin{enumerate}
\item $\theta^2 = \id$; and
\item $B_+^\theta = B_-$; moreover, every Borel subgroup $B'$ of $G$ containing $T = B_+ \cap N = B_- \cap N$ is mapped to an opposite Borel subgroup, i.e., $B' \cap \theta(B') = T$.
\end{enumerate}
Then $\theta$ centralizes the Weyl group $N/T$.
In particular, $\theta$ is a $BN$-flip.
\end{proposition}

\begin{proof}
For each $s \in S$ the set $P_s := B_+ \cup B_+sB_+$ is a rank one parabolic subgroup of positive sign of $G$. Let $n_s$ be a representative of $s$ in $N$. Then
$P_s^\theta = B_+^\theta \cup B_+^\theta n_s^\theta B_+^\theta
= B_- \cup B_- n_s^\theta B_-$
is a parabolic subgroup of negative sign of $G$, as it is the image of a subgroup of $G$ under the group automorphism $\theta$ and contains $B_-$. It consists of precisely two Bruhat double cosets, implying it must be a rank one parabolic subgroup. Hence $n^\theta_s$ is a representative of some $s'_{n_s} \in S$. As $s'$ is independent of the choice of $n_s$, the map $\theta$ induces a permutation of $S$ by setting $s^\theta:=s'_{n_s}$. Since $\theta$ maps every Borel group containing $T$ to an opposite one, for every $s \in S$ the chamber $sB_+$ is opposite to the chamber $s^\theta B_-$ in the associated twin building $\mc{B}$. That is, they have codistance $1_W$, whence $B_+ s^{-1} s^\theta B_- = B_+ B_-$ which by the uniqueness of the Birkhoff decomposition yields that $s^\theta=s$. Hence $\theta$ centralizes $W = \gen{S}$.
\end{proof}

\begin{corollary} \label{3.5}
Let $G$ be a group with a root group datum $\{ U_\alpha \}_{\alpha \in \Phi}$, let $B_+$, $B_-$, $N$ be the induced saturated twin $BN$-pair, and let $\theta$ be an automorphism of $G$ satisfying
\begin{enumerate}
\item $\theta^2 = \id$; and
\item $B_+^\theta = B_-$; moreover, every Borel subgroup $B'$ of $G$ containing $T = B_+ \cap N = B_- \cap N$ is mapped to an opposite Borel subgroup, i.e., $B' \cap \theta(B') = T$.
\end{enumerate}
Then $\theta$ centralizes the Weyl group $N/T$.
In particular, $\theta$ is a $BN$-flip. Moreover, $\theta$ normalizes $\langle U_\alpha, U_{-\alpha} \rangle$ for each (simple) root $\alpha$.
\end{corollary}

\begin{proof}
In view of the preceding proposition it remains to note that the final statement of the corollary follows from the main result of \cite{Caprace/Muehlherr:2006}, since $\theta$ maps bounded subgroups of $G$ to bounded subgroups of $G$, because it acts on the building.
\end{proof}

\begin{rem}
If in Corollary \ref{3.5} the root system $\Phi$ does not admit a direct summand of type $A_1$, then \cite[Lemma 8.17]{Abramenko/Brown:2008} also implies that $\theta$ preserves root subgroups and hence normalizes $\langle U_\alpha, U_{-\alpha} \rangle$ for each (simple) root $\alpha$, because $\theta$ acts on the building.
\end{rem}

\begin{proposition} \label{chevalley-is-flip}
Any $\sigma$-twisted Chevalley involution $\theta$ of a group $G$ is a $BN$-flip.
\end{proposition}

\begin{proof}
By definition, $\theta$ is an involution. Furthermore, the Borel subgroup $B_+$ is generated by $T$ and the set of root groups associated to the positive root system $\Phi_+\subset \Phi$. More precisely, $B_+=T.\gen{U_\alpha \mid \alpha \in \Phi_+}$. Since $T = \bigcap_{\alpha \in \Phi} N_G(U_\alpha)$ by \cite[Corollary 5.3]{Caprace/Remy:2008-LN}, the involution $\theta$ stabilizes $T$ and maps $B_+$ to $B_-=T.\gen{U_{-\alpha} \mid \alpha \in \Phi_+}$. Finally, $\theta$ acts trivially on $W = N/T$ as each root $\alpha$ of the root lattice of $W$ is mapped onto its negative $-\alpha$, which means that the reflection given by $\alpha$ is mapped onto the reflection given by $-\alpha$, which is identical to the reflection given by $\alpha$. 
\end{proof}

\begin{lemma} \label{chev-inv-existence}
Let $\F$ be a field with at least four elements, let $\sigma$ be an automorphism of $\F$ of order $1$ or $2$, let $G$ be a group with a $2$-spherical $\F$-locally split root group datum. Then $G$ admits a $\sigma$-twisted Chevalley involution.
\end{lemma}

\begin{proof}
By \cite{Abramenko/Muehlherr:1997} (and also by the unpublished manuscript \cite{Muehlherr}) the group $G$ is a universal enveloping group of the amalgam $\bigcup_{\alpha, \beta \in \Pi} X_{\alpha,\beta}$ for a system $\Pi$ of fundamental roots of $\Phi$. This means that any automorphism of $\bigcup_{\alpha, \beta \in \Pi} X_{\alpha,\beta}$ induces an automorphism of $G$.
For each pair $\alpha, \beta \in \Pi$ the $\sigma$-twisted Chevalley involution of the split semisimple algebraic group $X_{\alpha,\beta}$ induces $\sigma$-twisted Chevalley involutions $\theta_\alpha$ on $X_\alpha$ and $\theta_\beta$ on $X_\beta$. Therefore there exists an involution of the amalgam $\bigcup_{\alpha, \beta \in \Pi} X_{\alpha,\beta}$ inducing $\theta_\alpha$ on $X_\alpha$. Consequently there exists an involution $\theta$ on its universal enveloping group $G$ inducing $\theta_\alpha$ on each subgroup $X_\alpha$. This involution $\theta$ of $G$ by construction is a $\sigma$-twisted Chevalley involution of $G$.
\end{proof}

We now turn our attention to a local-to-global analysis of Iwasawa decompositions.
It is well-known that an adjacency-preserving action of a group $G$ on a connected chamber system $\mathcal{C}$ over $I$ is transitive if and only if there exists a chamber $c \in \mathcal{C}$ such that for each $i \in I $ the normalizer $N_G([c]_i)$ acts transitively on the $i$-panel $[c]_i$ of $\mathcal{C}$ containing $c$. It is implied by the following observation concerning permutation groups, whose proof is left to the reader. 

\begin{proposition}
Let $X$ be a set endowed with a family of equivalence relations $(\sim_i)_{i\in I} \subseteq X \times X$ such that the transitive hull of $\bigcup_{i \in I} \sim_i$ equals $X \times X$. Moreover, let $G$ be a group acting on $X$ as permutations preserving each equivalence relation. If there exists a point $p\in X$ such that the normalizer $G_i:=G_{[p]_i}$ of $[p]_i$ in $G$ acts transitively on $[p]_i$, then $G$ acts transitively on $X$.
\end{proposition}

\begin{corollary}
Let $\mathcal{C}$ be a connected chamber system and let $G$ be a group of automorphisms of $\mathcal{C}$. The group $G$ acts transitively on $\mathcal{C}$ if and only if there exists a chamber $c\in \mathcal{C}$, such that for each panel $P$ of $\mathcal{C}$ containing $c$ the stabilizer $G_P$ acts transitively on $P$.
\end{corollary}

\begin{corollary} \label{local-to-global}
Let $\mc{B}$ be a twin building obtained from a group $G$ with a twin $BN$-pair.
Let $\theta$ be a $BN$-flip of $G$, let $K:=\Gt$ be the group of all elements of $G$ centralized by $\theta$. The group $K$ acts transitively on the positive/negative half of $\mc{B}$ if and only if there exists a chamber $c$ opposite $c^\theta$ such that for each panel $P$ at $c$ the stabilizer $K_P$ of that panel in $K$ acts transitively on $P$.
\end{corollary}

We finally have assembled all tools required to prove our main result up to the rank one analysis conducted in the remainder of this article.

\begin{proof}[Proof of Theorem \ref{thm1}]
Assume the existence of an Iwasawa decomposition of $G$. By definition there exists
 an involution $\theta$ of $G$ such that $G = \Gt B_+$. Hence any Borel subgroup of $G$ is mapped onto an opposite one, so that by Corollary \ref{3.5} the involution $\theta$ centralizes the Weyl group $N/T$ and, for any simple root $\alpha$, normalizes the group $X_\alpha:=\gen{U_\alpha,U_{-\alpha}}$, which by $\F$-local splitness is isomorphic to $(\mathrm{P})\SL_2(\F)$.  In particular the restriction $\theta_{|X_\alpha}$ of $\theta$ to $X_\alpha$ is a $BN$-flip. 
 
We now argue that this restricted $BN$-flip induces an Iwasawa decomposition of $X_\alpha$. Let $P_\alpha$ be the panel of the building corresponding to the root $\alpha$. By Corollary \ref{local-to-global} we know that $(G_{P_\alpha})_\theta=G_{P_\alpha}\cap\Gt$ acts transitively on $P_\alpha$, and it remains to show that this is also the case for $(X_\alpha)_\theta=X_\alpha\cap\Gt$. 
First observe that $P_{-\alpha}=\theta(P_\alpha)$ and hence $(G_{P_\alpha})_\theta$ also stabilizes the panel $P_{-\alpha}$. For, if $g\in (G_{P_\alpha})_\theta$, then $g.P_{-\alpha} = g.\theta(P_\alpha)=\theta(g.P_\alpha)=\theta(P_\alpha)=P_{-\alpha}$ and so $g\in (G_{P_\alpha})_\theta=G_{P_\alpha}\cap G_{P_{-\alpha}}\cap\Gt$.
If $x\in(G_{P_\alpha})_\theta$ stabilizes the chamber $B_+$ in $P_\alpha$, then  $x.B_- = x.\theta(B_+)=\theta(x.B_+)=\theta(B_+)=B_-$.
We conclude that $x\in B_+\cap B_-=T$. Moreover, the group $U_\alpha<X_\alpha$ stabilizes $B_+$ and acts transitively on $P_{-\alpha}$. Thus, in fact $(G_{P_\alpha})_\theta= (X_\alpha T)\cap\Gt$.
Any $t\in T\setminus X_\alpha$ acts trivially on $P_\alpha$. Hence, since $(G_{P_\alpha})_\theta$ acts transitively on $P_\alpha$, so does $(X_\alpha)_\theta$. Accordingly $X_\alpha$ admits an Iwasawa decomposition.
 
Therefore, by Corollary \ref{sl2iwa} below, the field $\F$ admits an automorphism $\sigma$ with the required properties.

\medskip
For the converse implication, let $\theta$ be the $\sigma$-twisted Chevalley involution of $G$. For each $\alpha \in \Phi$ the involution $\theta$ induces a $BN$-flip $\theta_\alpha$ on $X_\alpha$.
By Proposition \ref{crit} below, these induced flips are transitive.  Hence by Corollary \ref{local-to-global}, we have $G = \Gt B_+$, proving that $G$ admits an Iwasawa decomposition.
\end{proof}

\section{Transitive involutions} \label{transinv}

As we have seen in the preceding section, the key to understanding Iwasawa decompositions lies in the structure of the rank one groups.
Since we are interested in groups with a locally split root group datum, we can restrict our attention to flip automorphisms of $\SL_2(\F)$ and $\PSL_2(\F)$ where $\F$ is an arbitrary field. 

In Section \ref{sect:classify-flips} we classify all suitable flip automorphisms of these two groups. We use that all automorphisms of $\PSL_2(\F)$ are induced via projection from automorphisms of $\SL_2(\F)$, which follows from the fact that $\SL_2$ is perfect, if $\abs{\F}\geq4$, and is easily verified over the fields of two and three elements. Alternatively, one can use the classification of endomorphisms of Steinberg groups or apply the results in \cite{Ren/Wan/Wu:1987}. Hence it suffices here to study flips of $\SL_2(\F)$. 

In Section \ref{sect:fixpointgroups} we compute the fixed point groups of these flips and give a geometric interpretation for a rank one Iwasawa decomposition. This finally enables us to give a nice sufficient and necessary algebraic criterion for such a local Iwasawa decomposition in Section \ref{sect:local-iwasawa}.

Furthermore, in Section \ref{sect:moufang}, we study flips of general Moufang sets which correspond to arbitrary rank one groups. Our aim is to show that it seems feasible to extend the theory to groups beyond $\F$-locally split ones. As a first step we present some results for $\SL_2(\mathbb{D})$ for a division ring $\mathbb{D}$.

\subsection{\boldmath Flip automorphisms of $G=\mathbf{SL_2(\F)}$} \label{sect:classify-flips}

In order to be able to understand flips of $G=\SL_2(\F)$ we need to specify a suitable root group datum of $G$. To this end consider $\mathrm{SL}_2(\F)$ as a matrix group acting on its natural module and let $T$ denote the subgroup of diagonal matrices, which is a maximal torus of $G$. Let $U_+$ and $U_-$ denote the subgroups of upper resp.\ lower triangular unipotent matrices, which are the root subgroups with respect to the root system of type $A_1$ associated to $T$. The standard Borel subgroups of $G$ then are the groups $B_+:=T.U_+$, $B_-:=T.U_-$. Finally, set $N:=N_G(T)$ to obtain a $BN$-pair.
Consider a $BN$-flip $\theta$ with respect to this $BN$-pair, i.e.\ an involutory automorphism $\theta$ of $G$ which interchanges $B_+$ and $B_-$ and centralizes $N/T$. It follows that $\theta$ stabilizes $T$ and interchanges $U_+$ and $U_-$. 

Let $K:=C_G(\theta)$, the fixed point group of $\theta$. Then $\theta$ induces an Iwasawa decomposition $G = K B_+$ if and only if  $K$ acts transitively on the projective line $\P_1(\F)=G/B_+$. In this case $\theta$ is called \Defn{transitive}.
Since $\theta$ interchanges $U_+$ and $U_-$ and since the root subgroups are isomorphic to $(\F,+)$, there must exist a group automorphism $\phi\in\Aut(\F,+)$ such that the equalities
\begin{equation} \label{eq:phi-induces-theta}
\theta\left(
\begin{pmatrix}
  1 & x \\
  0 & 1
\end{pmatrix}
\right)
=
\begin{pmatrix}
  1 & 0 \\
  \phi(x) & 1
\end{pmatrix}
\qquad\text{ and }\qquad
\theta\left(
\begin{pmatrix}
  1 & 0 \\
  y & 1
\end{pmatrix}
\right)
=
\begin{pmatrix}
  1 & \phi^{-1}(y) \\
  0 & 1
\end{pmatrix}
\end{equation}
hold. This weak assumption implies much stronger properties of $\phi$, as the next lemma shows.

\begin{lemma} \label{groupfield}
The group automorphism $\phi\in\Aut(\F,+)$ induces an involution $\theta\in\Aut(\SL_2(\F))$ if and only if $\phi(x)=\eps x^\sg$ for some field automorphism $\sigma\in\Aut(\F)$ of order $1$ or $2$ and some $\eps\in\Fix_\F(\sg)$.
\end{lemma}

\begin{proof}
The key here is to use the following equation derived from the Steinberg relations for Chevalley groups of rank one for various values $s,t\in\F^\times$, $u\in\F$:
\begin{equation} \label{steinbergRel}
\begin{pmatrix}
  1 & s \\
  0 & 1
\end{pmatrix}
\begin{pmatrix}
  1 & 0 \\
  -\frac{1}{s} & 1
\end{pmatrix}
\begin{pmatrix}
  1 & s-t+u \\
  0 & 1
\end{pmatrix}
\begin{pmatrix}
  1 & 0 \\
  \frac{1}{t} & 1
\end{pmatrix}
\begin{pmatrix}
  1 & -t \\
  0 & 1
\end{pmatrix}
=
\begin{pmatrix}
  \frac{s}{t} & 0 \\
  -\frac{u}{st} & \frac{t}{s}
\end{pmatrix}
.
\end{equation}
Assume $\theta$ is an involution induced by some group automorphism $\phi$ of $(\F,+)$ as described in (\ref{eq:phi-induces-theta}).
In case $s=t$ we can thus apply $\theta$ to (\ref{steinbergRel}) and obtain the equality
\begin{equation} \label{phiThetaRel}
\begin{pmatrix}
  1 & 0 \\
  \phi(t) & 1
\end{pmatrix}
\begin{pmatrix}
  1 & -\phi^{-1}(\frac{1}{t}) \\
  0 & 1
\end{pmatrix}
\begin{pmatrix}
  1 & 0 \\
  \phi(u) & 1
\end{pmatrix}
\begin{pmatrix}
  1 & \phi^{-1}(\frac{1}{t}) \\
  0 & 1
\end{pmatrix}
\begin{pmatrix}
  1 & 0 \\
  -\phi(t) & 1
\end{pmatrix}
=
\begin{pmatrix}
  1 & -\phi^{-1}(\frac{u}{t^2}) \\
  0 & 1
\end{pmatrix}
.
\end{equation}
Reducing (\ref{phiThetaRel}) further by setting $u=t$ and expanding the left side, we arrive at
\begin{equation*}
\begin{pmatrix}
  1-\phi(t)\phi^{-1}(\frac{1}{t}) \left(1- \phi(t)\phi^{-1}(\frac{1}{t}) \right)
         & -\phi(t)\phi^{-1}(\frac{1}{t})^2 \\
  \phi(t) \left(1-\phi(t)\phi^{-1}(\frac{1}{t})\right)^2
         & 1+\phi(t)\phi^{-1}(\frac{1}{t}) \left(1- \phi(t)\phi^{-1}(\frac{1}{t}) \right)
\end{pmatrix}
=
\begin{pmatrix}
  1 & - \phi^{-1}(\frac{1}{t}) \\
  0 & 1
\end{pmatrix}
\end{equation*}
which readily implies $\phi^{-1}(t)=\phi(\frac{1}{t})^{-1}$ for all $t\in\F^\times$. Defining $\eps:=\phi(1)$ and $\sg(x):=\frac{\phi(x)}{\eps}$ (note that $\eps\neq0$ since $\phi$ is an automorphism of the group $(\F,+)$), we again use (\ref{steinbergRel}), this time for $u=0$ and arbitrary $s$ and $t$. We obtain
\[
\theta
\begin{pmatrix}
  \frac{s}{t} & 0 \\
  0 & \frac{t}{s}
\end{pmatrix}
=
\begin{pmatrix}
  \frac{\phi(t)}{\phi(s)} & 0 \\
  0 & \frac{\phi(s)}{\phi(t)}
\end{pmatrix}
=
\begin{pmatrix}
  \frac{\sg(t)}{\sg(s)} & 0 \\
  0 & \frac{\sg(s)}{\sg(t)}
\end{pmatrix}
.
\]
We use this equality twice, once substituting $(y,xy)$ for $(s,t)$ and once $(1,x)$ for $(s,t)$ in order to obtain our final equality
\[
\begin{pmatrix}
  \frac{\sg(xy)}{\sg(y)} & 0 \\
  0 & \frac{\sg(y)}{\sg(xy)}
\end{pmatrix}
=
\theta
\begin{pmatrix}
  x^{-1} & 0 \\
  0 & x
\end{pmatrix}
=
\begin{pmatrix}
  \frac{\sg(x)}{\sg(1)} & 0 \\
  0 & \frac{\sg(1)}{\sg(x)}
\end{pmatrix}
=
\begin{pmatrix}
  \sg(x) & 0 \\
  0 & \frac{1}{\sg(x)}
\end{pmatrix}
.
\]
This allows us to conclude that $\sg(xy)=\sg(x)\sg(y)$. We already know that $\sg(x+y)=\sg(x)+\sg(y)$, $\sg(0)=0$, $\sg(1)=1$, and hence $\sg\in\Aut(\F)$ as required. Furthermore, $\sg(\eps)=\eps$ since
\[1=\sg(\eps \eps^{-1})
  =\sg(\eps)\sg(\eps^{-1})
  =\sg(\eps) \eps^{-1}\phi(\eps^{-1})
  = \eps^{-1}\sg(\eps) \phi^{-1}(\eps)^{-1}
  = \eps^{-1}\sg(\eps)
  .
 \]
Finally $\sg^2=\id$, since
\[\sg^{-1}(x)
 =\phi^{-1}(\eps x)
 = \frac{1}{\phi(\eps^{-1} x^{-1})}
 = \frac{1}{\eps \sg(\eps^{-1} x^{-1})}
 = \frac{1}{\sg(x^{-1})}
 = \sg(x)
 .
\]
The converse implication, deriving a group automorphism $\theta$ of $\SL_2(\F)$ from a given group automorphism $\phi$ of $(\F,+)$, results from the fact that the following automorphism restricts as required to $U_+$ resp.\ $U_-$:
\[ \theta : \SL_2(\F) \to \SL_2(\F) : X \mapsto \theta(X)=
\begin{pmatrix}
  0 & 1 \\
  \eps & 0
\end{pmatrix}
X^\sigma
\begin{pmatrix}
  0 & \eps^{-1} \\
  1 & 0
\end{pmatrix}.\tag*{\qedhere}
\]
\end{proof}
For an alternative proof see \cite[Section 3.1.1]{Horn:2008}.
\begin{defn}
For a field automorphism $\sigma$ of $\mathbb{F}$ of order $1$ or $2$ and $\delta \in \Fix_\F(\sigma)$ define 
\[ \theta_{\delta,\sigma} : \SL_2(\F) \to \SL_2(\F) : X \mapsto \theta_{\delta,\sigma}(X)=
\begin{pmatrix}
  0 & 1 \\
  -\delta^{-1} & 0
\end{pmatrix}
X^\sigma
\begin{pmatrix}
  0 & -\delta \\
  1 & 0
\end{pmatrix}.
\]
\end{defn}

By slight abuse of notation, we will use the same symbol $\theta_{\delta,\sigma}$ to denote the induced flip on $\PSL_2(\F)$.

\subsection{Centralizers of flips} \label{sect:fixpointgroups}

We now turn our attention to the centralizers of a given flip $\theta$. It is easy to verify that
\[ K_{\delta,\sg} :=C_{\SL_2(\F)}(\theta_{\delta,\sg})
  = \left\{
\begin{pmatrix}
  u^\sigma & \delta v^\sigma \\
  -v & u
\end{pmatrix}
\mid uu^\sigma+\delta vv^\sigma = 1
 \right\},\]
which is precisely the group preserving the $\sigma$-sesquilinear form
  \[ f(x,y):=x^T \begin{pmatrix} 1 & 0\\ 0 & \delta \end{pmatrix} y^\sigma \]
on the vector space $\F^2$ and its associated unitary form $q(x):=f(x,x)$.
This alternative characterization will turn out to be quite useful.

\medskip

For $\PSL_2(\F)$, the situation is slightly different. Let $Z$ denote the center of $\SL_2(\F)$. By definition $\PSL_2(\F)=\SL_2(\F)/Z$, so that the centralizer of $\theta$ in $\PSL_2(\F)$ is
$  C_{\PSL_2(\F)}(\theta) = \{gZ \in \PSL_2(\F) \mid (gZ)^\theta = gZ \}. $
We are mainly interested in the action of this centralizer on $\P_1(\F)$. Since the action of $\PSL_2(\F)$ is induced by that of $\SL_2(\F)$, this means studying the preimage of the centralizer in $\SL_2(\F)$. This suggests the following definition:
\begin{defn}
Let $\theta$ be an automorphism of $\SL_2(\F)$. We define the \Defn{projective centralizer} of $\theta$ in $\SL_2(\F)$ as the group
$PC_{\SL_2(\F)}(\theta) := \{g\in \SL_2(\F) \mid g^\theta \in gZ \}$,
which is the preimage of $C_{\PSL_2}(\theta)$ in $\SL_2(\F)$ under the canonical projection $\pi:\SL_2\to\PSL_2$.
\end{defn}

We compute
\[
  PK_{\delta,\sg} :=
  PC_{\SL_2(\F)}(\theta_{\delta,\sigma}) =
 \left\{
 \begin{pmatrix}
  \eps u^\sigma & \delta \eps v^\sigma \\
  -v & u
 \end{pmatrix}
 \mid uu^\sigma+\delta vv^\sigma = \eps,\; \eps \in \{-1,+1\}
 \right\}.
\]
While $K_{\delta,\sg}$ preserves the  $\sigma$-sesquilinear form $f(x,y)$ and its associated unitary form $q(x)$, the group $PK_{\delta,\sg}$ preserves these forms up to sign.

\subsection{Transitivity of centralizers of flips} \label{sect:local-iwasawa}

The observations of the previous sections allows us to characterize when $\theta_{\delta,\sg}$ is transitive. Set $N:\F\to\F:a\mapsto aa^\sigma$. Note that $q\vct{a}{b}=N(a)+\delta N(b)$. All results are written with $\PSL_2(\F)$ in mind. The corresponding results for $\SL_2(\F)$ can be obtained by replacing $PK_{\delta,\sg}$ by $K_{\delta,\sg}$ and substituting $1$ for $\eps$.

\begin{lemma} \label{conj}
If the involution $\theta_{\delta,\sg}$ is transitive, then it is conjugate to $\theta_{1,\sg}$ by an element of $\GL_2(\F)$ normalizing both $B_+$ and $B_-$. In case $\mathrm{char}(\F) = 2$, the involution $\theta_{\delta,\sg}$ cannot be transitive.
\end{lemma}

\begin{proof}
Due to the transitivity of $PK_{\delta,\sg}$ on $\P_1(\F)$, there exists
$g\in PK_{\delta,\sg}$ such that $g\vct{1}{0}=\vct{0}{x}$ for some $x\neq0$. Thus 
\[ 1=q\vct{1}{0}=\eps q(g\vct{1}{0})=\eps q\vct{0}{x}=
     \eps q\vct{0}{1}N(x)=\eps \delta N(x), \]
with $\eps\in\{-1,+1\}$, whence $\delta = \eps N(x^{-1})$.
If $\eps=-1$, then $q$ would be isotropic, as $q\vct{1}{x}=0$,
contradicting transitivity. (This in particular implies that in case $\mathrm{char}(\F) = 2$, the involution $\theta_{\delta,\sg}$ cannot be transitive.) Thus $\delta = N(x^{-1})$.
Let $Y:= \mtr{1}{0}{0}{x^{-1}}$,  and denote by $\mathrm{Inn}_Y(g)$
the inner automorphism of $G$ induced by $Y$. Then conjugating
$\theta_{1,\sg}$ by $\mathrm{Inn}_Y(g)$ yields $\theta_{\delta,\sg}$.
\end{proof}

Because of the preceding lemma it remains to determine when exactly $\theta_{1,\sg}$ is transitive.

\begin{prop} \label{crit}
The involution $\theta_{1,\sg}$ is transitive if and only if $-1$ is not a norm, and the sum of two norms is $\eps$ times a norm, where $\eps\in\{+1,-1\}$.
\end{prop}

\begin{proof}
Assume $\theta_{1,\sg}$ is transitive. Take an arbitrary nonzero vector
$\vct{a}{b}\in\F^2$. Due to transitivity, there exists $g\in PK_{1,\sg}$
such that $g\vct{a}{b}=\vct{x}{0}$ for some nonzero $x\in\F$. Consequently,
\[N(a)+N(b)=q\vct{a}{b}=\eps q(g\vct{a}{b})=\eps q\vct{x}{0}=\eps N(x)\in N(\F),
  \; \eps \in \{-1,+1\}, \]
proving that a sum of two norms is a norm or $-1$ times a norm. Furthermore, $-1$ is not a norm as else there would be $x\in\F$ with $N(x)=-1$ and we would obtain the isotropic vector $\vct{1}{x}$ contradicting transitivity.

For the converse implication assume that the sum of two norms is a norm or $-1$ times a norm, and $-1$ is not a norm. It suffices to show that there exists an element of $PK_{1,\sg}$ that maps an arbitrary non-trivial vector $\vct{a}{b} \in \F^2$ onto some non-trivial vector $\vct{x}{0}$. Choose $x$ such that $\eps N(x)=N(a)+N(b)$ for $\eps\in\{-1,+1\}$, and note that $x\neq0$, as else $N(a)=-N(b)$ contradicting that $-1$ is not a norm.
Thus the equation 
\[
   \begin{pmatrix} \eps\left(\frac{a}{x}\right)^\sigma & \eps\left(\frac{b}{x}\right)^\sigma \\ -\frac{b}{x} & \frac{a}{x} \end{pmatrix}
  \begin{pmatrix} a \\ b \end{pmatrix}
  = \begin{pmatrix} x \\ 0 \end{pmatrix}
\]
finishes the proof, since the given matrix is clearly in $PK_{1,\sg}$.
\end{proof}

\begin{cor} \label{sl2iwa}
The group $(\mathrm{P})\SL_2(\F)$ admits an Iwasawa decomposition if and only if $\F$ admits an automorphism $\sigma$ of order $1$ or $2$ such that 
\begin{enumerate}
\item $-1$ is not a norm, and
\item
  \begin{enumerate}
    \item either a sum of norms is a norm (in the $\SL_2(\F)$ case)
    \item or a sum of norms is $\eps$ times a norm, where $\eps\in\{+1,-1\}$
          (in the $\PSL_2(\F)$ case),
  \end{enumerate}
\end{enumerate}
 with respect to the norm map $N : \F \to \Fix_\F(\sigma) : x \mapsto x x^\sigma$.
\end{cor}
\begin{proof}
Assume we have an Iwasawa decomposition of $G$. Then we have an involution $\theta$ which interchanges $U_+$ and $U_-$ and satisfies $G=\Gt B_+$, whence $\theta$ is transitive. Then by Lemmas \ref{groupfield} and \ref{conj} plus Proposition \ref{crit} the claim for $\F$ follows. If on the other hand $\F$ is as described, then again by Proposition \ref{crit} the map $\theta_{1,\sigma}$ induces an Iwasawa decomposition.
\end{proof}

\subsection{Fields permitting Iwasawa decompositions}

Besides the real closed fields and the field of complex numbers there exist lots of fields admitting automorphisms that satisfy the conditions of Corollary \ref{sl2iwa}. For instance any pythagorean formally real field $\F$ satisfies the conditions of Corollary \ref{sl2iwa} with respect to the identity automorphism as does $\F[\sqrt{-1}]$ with respect to the non-trivial Galois automorphism. Such fields have been studied very thoroughly, cf.\ \cite{Lam:1973}, \cite{Lam:2005}, \cite{Rijwade:1993}. In the $\mathrm{PSL}_2(\F)$ case, the finite fields $\mathbb{F}_q$ with $q \equiv 3 \mod 4$ yield additional examples.

Inspired by classical Lie theory and the passage from complex Lie groups to their split real forms, the question arises when an Iwasawa decomposition $G = \Gt B_+$ of a group $G$ with an $\F$-locally split root group datum with respect to an involution $\theta$ involving a non-trivial field automorphism $\sigma : \F \to \F$ implies the existence of an Iwasawa decomposition over the field $\Fix_\F(\sigma)$ with respect to an involution involving the trivial field automorphism on $\Fix_\F(\sigma)$. The following example shows that this in general is not the case.

\begin{ex}
Let $\F$ be a formally real field which is not pythagorean and admits four square classes. Such fields exist, see for example \cite{Szymiczek:1975}. This means exactly two square classes contain absolutely positive elements, so that there exists a unique ordering. 
Choose a positive non-square element $w\in\F$. Set $\alpha:=\sqrt{-w}$ and $\tilde\F:=\F[\alpha]$. Then
\[ N(x_0+\alpha x_1) + N(y_0+\alpha y_1) = x_0^2+wx_1^2+y_0^2+wy_1^2\]
which is a non-negative number, hence either a square or a square multiple of $w$. Hence there exist $z_0$ and $z_1$ in $\F$ such that
\[ N(x_0+\alpha x_1) + N(y_0+\alpha y_1) = x_0^2+wx_1^2+y_0^2+wy_1^2
= z_0^2+wz_1^2 = N(z_0+\alpha z_1)\]
and thus the field $\tilde\F$ together with the non-trivial Galois automorphism satisfies the conditions of Corollary \ref{sl2iwa}, while $\F$ together with the identity does not, because $\F$ is not pythagorean.
\end{ex}

\section{Moufang sets} \label{sect:moufang}

In this last section, we turn our attention to the study of Moufang sets.
Our motivation is Corollary \ref{3.5} which states that in order to understand flips of arbitrary non-split groups with a root group datum one needs to understand flips of arbitrary rank one subgroups. Moufang sets are essentially equivalent to these and seem to be the natural setting to study flips. 
For a concise introduction to Moufang sets, we refer the interested reader to \cite{DS09}.

Maybe Moufang sets, or rather Moufang sets with transitive flips (cf.\ Sections \ref{sect:moufang-flips} and \ref{sect:moufang-trans}), can be used for a better understanding of certain anisotropic forms of reductive algebraic groups. Indeed, the machinery developed in this paper implies that an anisotropic form $K$ of a reductive algebraic group $G$ acts transitively on the building $\mathcal{B}$ of $G$ if and only if there exists a chamber $c$ of $\mathcal{B}$ such that, for each $i$-panel $P$ containing $c$, the stabilizer $K_P$ acts transitively on $P$, cf.\@ Corollary~\ref{local-to-global}. Therefore any result about Moufang sets with transitive flips has immediate consequences for certain classes of anisotropic algebraic groups, namely those which act transitively on a building.

\subsection{Moufang sets and pointed Moufang sets} \label{sect:moufang-general}

In order to be consistent with the standard notation used in the theory of Moufang sets we will always denote the action of a permutation on a set on the right, i.e.\@ we will write
$a\varphi$ rather than $\varphi(a)$.

\begin{defn}
        A \Defn{Moufang set} is a set $X$ together with a collection of subgroups $(U_x)_{x \in X}$,
        such that each $U_x$ is a subgroup of $\Sym(X)$ fixing $x$ and acting regularly
        (i.e.\@ sharply transitively) on $X \setminus \{ x \}$, and such that each $U_x$
        permutes the set $\{ U_y \mid y \in X \}$ by conjugation.
        The group $G := \langle U_x \mid x \in X \rangle$ is called the
        \Defn{little projective group} of the Moufang set;
        the groups $U_x$ are called \Defn{root groups}.
\end{defn}

Our approach to Moufang sets is taken from \cite{DW}.
Let $\mouf = (X, (U_x)_{x \in X})$ be an arbitrary Moufang set, and assume that two of the elements
of $X$ are called $0$ and $\infty$.
Let $U := X \setminus \{ \infty \}$.
Each $\alpha \in U_\infty$ is uniquely determined by the image of $0$ under $\alpha$.
If $0\alpha = a$, we write $\alpha =: \alpha_a$.
Hence $U_\infty = \{ \alpha_a \mid a \in U \}$.
We make $U$ into a (not necessarily abelian) group with composition $+$ and identity 0,
by setting
\begin{equation}\label{eq:alpha}
        a + b := a\alpha_b \,.
\end{equation}
Clearly, $U \cong U_\infty$.
Now let $\tau$ be an element of $G$ interchanging $0$ and $\infty$.
(Such an element always exists, since $G$ is doubly transitive on $X$.)
By the definition of a Moufang set, we have
\begin{equation}\label{eq:defU}
        U_0 = U_\infty^\tau \text{ and } U_a = U_0^{\alpha_a}
\end{equation}
for all $a\in U$.
In particular, the Moufang set $\mouf$ is completely determined by the group $U$ and the permutation $\tau$;
we will denote it by $\mouf = \mouf(U, \tau)$.

\begin{rem}\label{re:infty}
        In view of equation \eqref{eq:alpha}, it makes sense to use the convention that 
        $a + \infty = \infty + a = \infty$ for all $a \in U$.
\end{rem}

\begin{defn}\label{def:gamma}
        For each $a \in U$, we define $\gamma_a := \alpha_a^\tau$,
        i.e.\@ $x \gamma_a = (x \tau^{-1} + a) \tau$ for all $x \in X$.
        Consequently, $U_0 = \{ \gamma_a \mid a \in U \}$.
\end{defn}

The following two definitions may appear technical at first sight, but it turns out that they play a key role in the theory of Moufang sets.

\begin{defn}
        For each $a \in U^* = U \backslash \{ 0 \}$, we define a \Defn{Hua map} to be $$h_a := \tau \alpha_a \tau^{-1} \alpha_{-(a\tau^{-1})} \tau \alpha_{-(-(a\tau^{-1}))\tau} \in \Sym(X);$$
        if we use the convention of Remark~\ref{re:infty}, then we can
        write this explicitly as $h_a : X \to X : x \mapsto \bigl( (x\tau + a)\tau^{-1} - a\tau^{-1} \bigr) \tau -
                \bigl( -(a\tau^{-1}) \bigr) \tau$.
        Observe that each $h_a$ fixes the elements $0$ and $\infty$.
        We define the \Defn{Hua subgroup} of $\mouf$ as $H := \langle h_a \mid a \in U^* \rangle$.
        By \cite[Theorem~3.1]{DW}, the group $H$ equals $G_{0, \infty}:=\mathrm{Stab}_G(0,\infty)$, and by \cite[Theorem~3.2]{DW},
        the restriction of each Hua map to $U$ is additive, i.e.\@ $H \leq \Aut(U)$.
\end{defn}

\begin{defn}
        For each $a \in U^*$, we define a \Defn{$\mu$-map}
	$\mu_a:=\gamma_{(-a)\tau^{-1}}\alpha_a\gamma_{a\tau^{-1}}^{-1}$\,.
\end{defn}
\begin{lemma}
	For each $a \in U^*$, we have
        \begin{itemize}
            \item[\rm (i)]
		$\mu_a$ is the unique element in the set $U_0^* \alpha_a U_0^*$ interchanging $0$ and $\infty$\,;
            \item[\rm (ii)]
		$\mu_a^{-1} = \mu_{-a}$\,;
            \item[\rm (iii)]
		$\mu_a = \tau^{-1} h_a$\,.
        \end{itemize}
\end{lemma}
\begin{proof}
	See \cite[Section 3]{DS}.
\end{proof}

\begin{defn}\label{def:iso}
        Let $(X, (U_x)_{x \in X})$ and $(Y, (V_y)_{y \in Y})$ be two Moufang sets.
        A bijection $\beta$ from $X$ to $Y$ is called an \Defn{isomorphism} of Moufang sets,
        if the induced map $\chi_\beta : \Sym(X) \to \Sym(Y) : g \mapsto \beta^{-1} g \beta$
        maps each root group $U_x$ isomorphically onto the corresponding root group $V_{x\beta}$.
        An \Defn{automorphism} of $\mouf = (X, (U_x)_{x \in X})$ is an isomorphism from
        $\mouf$ to itself.
        The group of all automorphisms of $\mouf$ will be denoted by $\Aut(\mouf)$.
\end{defn}

Now we introduce pointed Moufang sets, which will be Moufang sets with a fixed identity element.
We will then, in analogy with the theory of Jordan algebras, introduce the notions of an isotope of
a pointed Moufang set, and we will define Jordan isomorphisms between Moufang sets.

\begin{defn}
        A \Defn{pointed Moufang set} is a pair $(\mouf, e)$, where $\mouf = \mouf(U, \rho)$
        is a Moufang set and $e$ is an arbitrary element of $U^*$.
        The \Defn{$\tau$-map} of this pointed Moufang set is $\tau := \mu_{-e} = \mu_e^{-1}$,
        and the \Defn{Hua maps} are the maps $h_a := \tau \mu_a = \mu_{-e} \mu_a$ for all $a \in U^*$.
        We also define the \Defn{opposite Hua maps} $g_a := \tau^{-1} \mu_a = \mu_e \mu_a$ for all $a \in U^*$.
        Clearly, $\mouf = \mouf(U, \tau) = \mouf(U, \tau^{-1})$.
\end{defn}

        Note that, in contrast with Moufang sets which are not pointed, the maps $\tau$, $h_a$ and $g_a$
        are completely determined by the data $(\mouf, e)$.
        On the other hand, there can be many different elements $f$ for which $(\mouf, e) = (\mouf, f)$,
        namely all those for which $\mu_e = \mu_f$.

\begin{defn}
        Let $(\mouf, e)$ and $(\mouf', f)$ be two pointed Moufang sets, with $\mouf = \mouf(U, \rho)$
        and $\mouf' = \mouf(U', \rho')$.
        A \Defn{pointed isomorphism} from $(\mouf, e)$ to $(\mouf', f)$
        is an isomorphism from $U$ to $U'$ mapping $e$ to $f$ and extending to a Moufang set isomorphism
        from $\mouf$ to $\mouf'$ (by mapping $\infty$ to $\infty'$).
        A pointed isomorphism from $(\mouf, e)$ to itself is called a \Defn{pointed automorphism}
        of $(\mouf, e)$, and the group of all pointed automorphisms is denoted by $\Aut(\mouf, e)$.
\end{defn}

        Observe that $G \cap \Aut(\mouf, e) = C_H(e)$.

\begin{defn}
        Let $(\mouf, e)$ be a pointed Moufang set, and let $a \in U^*$ be arbitrary.
        Then $(\mouf, a)$ is called the \Defn{$a$-isotope} of $(\mouf, e)$, or simply an \Defn{isotope}
        if one does not want to specify the element $a$.
        The $\tau$-map and the Hua maps of $(\mouf, a)$ will be denoted by
        $\tau^{(a)}$ and $h_b^{(a)}$, respectively.
        Observe that
        \begin{equation}\label{eq:isot}
                \tau^{(a)} = \mu_{-a} \text{ \ and \ } h_b^{(a)} = \mu_{-a} \mu_b = h_a^{-1} h_b
        \end{equation}
        for all $a,b \in U^*$.
\end{defn}

        Our notion of an $a$-isotope is, in a certain sense, the inverse of the usual notion of an $a$-isotope
        in (quadratic) Jordan algebras, where our $a$-isotope would be called the $a^{-1}$-isotope
        (where $a^{-1}$ denotes the inverse in the Jordan algebra) and where $h_b^{(a)} := h_a h_b$.
        It is, in the general context of Moufang sets, not natural to try to be compatible with this convention,
        because $h_a^{-1}$ is in general not of the form $h_b$ for some $b \in U^*$.
        In fact, we have $h_a^{-1} = g_{a\tau}$ for all $a \in U^*$;
        see \cite[Lemma~3.8(i)]{DW}.


\begin{defn}
        Let $(\mouf, e)$ and $(\mouf', f)$ be two pointed Moufang sets with $\mouf = \mouf(U, \rho)$
        and $\mouf' = \mouf(U', \rho')$, and with Hua maps $h_a$ and $k_a$, respectively.
        An isomorphism $\varphi$ from $U$ to $U'$ is called a \Defn{Jordan isomorphism} if
        $(bh_a)\varphi = (b\varphi)k_{a\varphi}$
        for all $a,b \in U^*$.
        If $(\mouf', f)$ is an isotope $(\mouf, a)$ of $(\mouf, e)$, then a Jordan isomorphism
        from $(\mouf, e)$ to $(\mouf, a)$ is called an \Defn{isotopy} from $(\mouf, e)$ to its $a$-isotope.
        Explicitly, a map $\varphi \in \Aut(U)$ is an isotopy if and only if
        \begin{equation}
                h_a \varphi = \varphi h_{a\varphi}^{(e\varphi)}
        \end{equation}
        for all $a \in U^*$.
        The group of all isotopies from $(\mouf, e)$ to an isotope is called the
        \Defn{structure group} of $(\mouf, e)$, and is denoted by $\Str(\mouf, e)$.
        Note that it is not clear whether $\Str(\mouf, e) \leq \Aut(\mouf)$.
        Also observe that $G \cap \Str(\mouf, e) = H$;
        we call $H$ the \Defn{inner structure group} of $(\mouf, e)$.
\end{defn}

\subsection{Flips of Moufang sets} \label{sect:moufang-flips}

Our goal in this section is to determine all involutions $\theta \in \Aut(G)$ interchanging $U_\infty$ and $U_0$.
Such an involution $\theta$ maps each $\alpha_a$ to some $\gamma_{a\varphi}$
and each $\gamma_b$ to some $\alpha_{b\psi}$. Since $\theta \in \Aut(G)$, we have $\varphi,\psi \in \Aut(U)$.
Moreover, $\theta^2 = \mathrm{id}$ implies $\psi = \varphi^{-1}$.
In particular, $\theta$ is competely determined by $\varphi$.
More precisely, for each $\varphi \in \Aut(U)$, we define
\[ \theta_\varphi : U_\infty \cup U_0 \to U_0 \cup U_\infty : \begin{cases}
        \alpha_a \mapsto \gamma_{a\varphi} \\
        \gamma_a \mapsto \alpha_{a\varphi^{-1}}
\end{cases} ; \]
the question is when $\theta_\varphi$ extends to an automorphism of $G$.
Observe that if $\theta_\varphi$ extends, then this extension is unique and is involutory,
since $\theta$ is involutory on $U_\infty \cup U_0$ and $G = \langle U_\infty \cup U_0 \rangle$.

\begin{proposition}\label{pr:phitau}
        Let $\varphi \in \Aut(U)$.
        Then $\theta_\varphi$ extends to an (involutory) automorphism of $G$ if and only if $(\varphi\tau)^2 = \mathrm{id}$.
        Moreover, if this is the case, then $\varphi \in \Aut(\mouf)$.
\end{proposition}
\begin{proof}
        Let $\theta := \theta_\varphi$ and $\beta := \varphi\tau$.
        Assume first that $\theta$ extends to an automorphism $\chi$ of $G$.
        Then
        \begin{equation}\label{eq:Ua}
                \chi(U_a) = \chi(U_0^{\alpha_a}) = \chi(U_0)^{\chi(\alpha_a)}
                = U_\infty^{\gamma_{a\varphi}} = U_{a\varphi\tau} = U_{a\beta}
        \end{equation}
        for all $a \in U$.
        Since $\theta^2$ is the identity on $U_\infty \cup U_0$ and since $G = \langle U_\infty, U_0 \rangle$,
        this implies that $\chi^2 = 1$ and hence $\beta^2 = 1$.

        Conversely, assume that $\beta^2 = 1$, and let $\chi_\beta$ be as in Definition~\ref{def:iso}.
        Then for all $a \in U$,
        \begin{align*}
                \chi_\beta(\alpha_a) 
                &= \alpha_a^{\varphi\tau} = \alpha_{a\varphi}^\tau = \gamma_{a\varphi} \,, \\
                \chi_\beta(\gamma_a)
                &= \gamma_a^{\varphi\tau} = \gamma_a^{\tau^{-1} \varphi^{-1}} = \alpha_a^{\varphi^{-1}}
                = \alpha_{a\varphi^{-1}} \,;
        \end{align*}
        hence $\chi_\beta$ and $\theta$ coincide on $U_\infty \cup U_0$.
        Note that $\chi_\beta$ is an (inner) automorphism of $\Sym(X)$, and hence
        the same calculation as in equation~\eqref{eq:Ua} (with $\chi_\beta$ in place of $\chi$)
        shows that $\beta \in \Aut(\mouf)$.
        Hence the restriction of $\chi_\beta$ to $G$ is an automorphism of $G$;
        this is the (unique) extension of $\theta$ to an element of $\Aut(G)$.

        Finally, since we have just shown that $\beta \in \Aut(\mouf)$ and
        since obviously $\tau \in \Aut(\mouf)$, we conclude that $\varphi \in \Aut(\mouf)$ as well.
\end{proof}

\begin{defn}
        An automorphism $\varphi \in \Aut(U)$ with the property that $(\varphi\tau)^2 = 1$ will be
        called a \Defn{flip automorphism} of $\mouf$.
\end{defn}

The following theorem gives important information about such flip automorphisms.

\begin{thm}\label{th:main}
        Let $\mouf$ be a Moufang set, and let $\varphi$ be a flip automorphism of $\mouf$.
        Then
        \[ g_{a\varphi} = \varphi \cdot h_a \cdot \varphi \]
        for all $a \in U^*$.
        Moreover, if $e$ is an identity element of $\mouf$, i.e.\@ $\tau = \mu_{-e}$,
        then $\varphi \in \Str(\mouf, e) \cap \Aut(\mouf)$.
\end{thm}
\begin{proof}
        For each $a \in U^*$, the map $g_a$ is the Hua map of $a$ with $\tau$ replaced by $\tau^{-1}$, and hence
        $g_{a\varphi} = \tau^{-1} \alpha_{a\varphi} \tau \alpha_{-a\varphi\tau}
                \tau^{-1} \alpha_{-(-(-a\varphi\tau))\tau^{-1}}$
        for all $a \in U^*$.
        Using the facts that $\alpha_a^\varphi = \alpha_{a\varphi}$, $\varphi\tau = \tau^{-1} \varphi^{-1}$
        and $(-a)\varphi = -a\varphi$ several times, we get
        $\varphi^{-1} g_{a\varphi} = \tau \alpha_{a} \tau^{-1} \alpha_{-a\tau^{-1}}
                \tau \alpha_{-(-(-a\tau^{-1}))\tau} \varphi = h_a \varphi$.
        In particular, if $e$ is an identity element of $\mouf$, then $h_e = 1$
        and hence $\varphi^{-1} g_{e\varphi} = \varphi$.
        It follows that $\varphi g_{e\varphi}^{-1} g_{a\varphi} = h_a \varphi$ for all $a \in U^*$.
        However,
        $g_{e\varphi}^{-1} g_{a\varphi} = (\mu_e \mu_{e\varphi})^{-1} (\mu_e \mu_{a\varphi})
                = (\mu_{-e} \mu_{e\varphi})^{-1} (\mu_{-e} \mu_{a\varphi}) = h_{e\varphi}^{-1} h_{a\varphi}
                = h_{a\varphi}^{(e\varphi)}$
        and hence $h_a \varphi = \varphi h_{a\varphi}^{(e\varphi)}$ for all $a \in U^*$,
        proving that $\varphi \in \Str(\mouf, e)$.
        The fact that $\varphi \in \Aut(\mouf)$ was shown in Proposition~\ref{pr:phitau} above.
\end{proof}


We will now illustrate the strength of Theorem~\ref{th:main} by explicitly determining all flips of $\PSL_2(\mathbb{D})$,
where $\mathbb{D}$ is a field or a skew field. This can be considered as a natural extension of the results from Lemma \ref{groupfield} to the non-commutative case.

\begin{proposition}\label{pr:D}
        Let $\mathbb{D}$ be an arbitrary field or skew field, and let $\mouf = \mouf(\mathbb{D})$ be the corresponding Moufang set,
        i.e.\@ the Moufang set $\mouf = \mouf(U, \tau)$ where $U := (\mathbb{D}, +)$ and
        $\tau : \mathbb{D}^* \to \mathbb{D}^* : x \mapsto -x^{-1}$.
        \begin{itemize}
            \item[\rm (i)]
                Let $\varphi$ be a flip automorphism of $\mouf$.
                Then there exists an automorphism or anti-automorphism $\sigma$ of $\mathbb{D}$ and an element
                $\eps \in \Fix_\mathbb{D}(\sigma)$ such that $x\varphi = \eps \sigma(x)$ for all $x \in \mathbb{D}$.
                If $\sigma$ is an automorphism, then $\sigma^2(x) = \eps^{-1} x \eps$ for all $x \in \mathbb{D}$;
                if $\sigma$ is an anti-automorphism, then $\sigma^2 = 1$.
            \item[\rm (ii)]
                Conversely, suppose that either $\sigma$ is an anti-automorphism of order~$2$ and
                $\eps \in \Fix_\mathbb{D}(\sigma)$ is arbitrary, or $\sigma$ is an automorphism
                such that $\sigma^2(x) = \eps^{-1} x \eps$ for some $\eps \in \Fix_\mathbb{D}(\sigma)$.
                Then the map $\varphi : \mathbb{D} \to \mathbb{D} : x \mapsto \eps \sigma(x)$ is a flip automorphism of $\mouf$.
        \end{itemize}
\end{proposition}
\begin{proof}
        \begin{compactenum}\itemsep1ex
            \item[(i)]
                Observe that $1 \in \mathbb{D}^*$ is an identity element of $\mouf$; also note that $\tau^2 = \mathrm{id}$.
                For all $a,b \in U^*$, we have $bh_a = aba$.
                The condition $(\varphi\tau)^2 = 1$ translates to
                \begin{equation}\label{eq:invD}
                        (a^{-1})\varphi = (a\varphi^{-1})^{-1}
                \end{equation}
                for all $a \in \mathbb{D}^*$.
                Let $\eps := 1\varphi$;
                then $bh_a^{(1\varphi)} = b h_{1\varphi}^{-1} h_a = a \eps^{-1} b \eps^{-1} a$
                for all $a,b \in U^*$.
                By Theorem~\ref{th:main}, $\varphi \in \Str(\mouf, e)$, which means that
                $bh_a \varphi = b\varphi h_{a\varphi}^{(1\varphi)}$
                for all $a,b \in U^*$, or explicitly,
                $(aba) \varphi = a\varphi \cdot \eps^{-1} \cdot b\varphi \cdot \eps^{-1} \cdot a\varphi$
                for all $a,b \in \mathbb{D}^*$.
                Now let $\sigma(a) := \eps^{-1} \cdot a\varphi$ for all $a \in \mathbb{D}$. Then $\sigma \in \Aut(U)$,
                and the previous equation can be rewritten as
                $\sigma(aba) = \sigma(a)\sigma(b)\sigma(a)$
                for all $a,b \in \mathbb{D}$, i.e.\@ $\sigma$ is a Jordan automorphism of $\mathbb{D}$.
                It is a well known result by Jacobson and Rickart \cite{Jacobson/Rickart:1950} (see also \cite[page 2]{Jacobson:1968}), which simply amounts to calculating that
                $\bigl( \sigma(ab) - \sigma(a) \sigma(b) \bigr) \cdot
                        \bigl( \sigma(ab) - \sigma(b)\sigma(a) \bigr) = 0$,
                that $\sigma$ is either an automorphism or an anti-automorphism of $\mathbb{D}$.
                Now by equation~\eqref{eq:invD}, we have
                $(\eps^{-1})\varphi = (\eps\varphi^{-1})^{-1} = 1^{-1} = 1$,
                and hence $\sigma(\eps^{-1}) = \eps^{-1}$; since $\sigma$ is an automorphism or
                anti-\hspace{0pt}automorphism,
                it follows that $\sigma(\eps) = \eps$.
                Finally, again by equation~\eqref{eq:invD}, we obtain
                $\sigma(\eps \sigma(a)) = \sigma(a\varphi) = \sigma \bigl( ((a^{-1})\varphi^{-1})^{-1} \bigr) 
                                = \sigma((a^{-1})\varphi^{-1})^{-1} = (\eps^{-1} a^{-1})^{-1} = a \eps$
                for all $a \in \mathbb{D}^*$.
                If $\sigma$ is an automorphism, then this can be rewritten as $\eps \sigma^2(a) = a \eps$,
                or $\sigma^2(a) = \eps^{-1} a \eps$;
                if $\sigma$ is an anti-automorphism, we get $\sigma^2(a) \eps = a \eps$, i.e.\@ $\sigma^2 = 1$.
            \item[(ii)]
                It suffices to check that equation~\eqref{eq:invD} holds.
                This amounts to checking that $\eps \sigma(a^{-1}) = (\sigma^{-1}(\eps^{-1} a))^{-1}$
                for all $a \in \mathbb{D}$.
                It is straightforward to check that this is valid in both cases.
        \qedhere
        \end{compactenum}
\end{proof}

By \cite{Ren/Wan/Wu:1987} the flips of $\mathrm{SL}_2(\mathbb{D})$ are just the lifts of the flips of $\mathrm{PSL}_2(\mathbb{D})$.

\subsection{Transitivity of the obvious flip} \label{sect:moufang-trans}

\begin{defn}
        If $\tau^2=\mathrm{id}$, then $\varphi=1$ is a flip automorphism.
        We will call the corresponding automorphism $\theta_1$ of $G$ (as defined in the beginning of this section)
        the \Defn{obvious flip}.
        Observe that $\theta_1$ is just conjugation by $\tau$.
\end{defn}

\begin{defn}
        A flip automorphism $\varphi \in \Aut(U)$ is called \Defn{transitive} if the group $\Fix_G(\theta_\varphi)$
        is transitive on $X$.
\end{defn}

        Let $\mouf = \mouf(U, \tau)$ be a Moufang set with $\tau^2 = \mathrm{id}$.
        Then the obvious flip $\theta_1$ is transitive if and only if
        $C_G(\tau)$ is transitive on $X$, because $\Fix_G(\theta_1) = C_G(\tau)$.

\begin{lemma}
        Let $\mouf = \mouf(U, \tau)$ be a Moufang set with $\tau^2 = \mathrm{id}$, and
        assume that the obvious flip is transitive.
        Then $\tau$ has no fixed points.
\end{lemma}
\begin{proof}
        Assume that $a\tau = a$ for some $a \in U^*$.
        Let $g \in C_G(\tau)$ be such that $0g = a$.
        Then $\infty g = 0\tau g = 0 g \tau = a\tau = a = 0g$ and hence $\infty = 0$, a contradiction.
\end{proof}

We now examine the transitivity of the obvious flip for $\mouf(\mathbb{D})$ where $\mathbb{D}$ is an arbitrary skew field.

\begin{defn}
        If $g = \mtr{a}{b}{c}{d} \in \GL_2(\mathbb{D})$, then
        the Dieudonn\'e determinant $\det(g) \in \mathbb{D}^* / [\mathbb{D}^*, \mathbb{D}^*]$ is defined as
        \[ \det(g) := \begin{vmatrix} a & b \\ c & d \end{vmatrix} := \begin{cases}
                ad - aca^{-1}b & \text{ if $a \neq 0$} \,; \\
                - cb & \text{ if $a = 0$} \,;
        \end{cases} \]
see \cite{Dieudonne:1943}.
        Then $\SL_2(\mathbb{D})$ is precisely the kernel of the Dieudonn\'e determinant, i.e.\@ 
        a matrix $g \in \GL_2(\mathbb{D})$ lies in $\SL_2(\mathbb{D})$ if and only if $\det(g) \in [\mathbb{D}^*, \mathbb{D}^*]$.
        Also observe that $\det(\lambda g) \equiv \det(g \lambda) \equiv \lambda^2 \det(g) \mod{[\mathbb{D}^*, \mathbb{D}^*]}$ for all $\lambda \in \mathbb{D}^*$.
\end{defn}

\begin{lemma}\label{le:DCG}
        Let $G = \SL_2(\mathbb{D})$ and let $\tau = \mtr{0}{1}{-1}{0} \in G$.
        Then
        \begin{align*}
                C_G(\tau) &= \left\{ \begin{pmatrix} a & b \\ -b & a \end{pmatrix} \Big\arrowvert
                        \begin{array}{lr} a^2 + aba^{-1}b \in [\mathbb{D}^*, \mathbb{D}^*] & \mbox{if } a \neq 0 \\
                        b^2 \in [\mathbb{D}^*, \mathbb{D}^*] & \mbox{if } a = 0 \end{array}\right\} \,; \\
                PC_G(\tau) &= \left\{ \begin{pmatrix} \epsilon a & \epsilon b \\ -b & a \end{pmatrix} \Big\arrowvert
                        \begin{array}{lr} \epsilon \cdot (a^2 + aba^{-1}b) \in [\mathbb{D}^*, \mathbb{D}^*] & \mbox{if } a \neq 0 \\
                        \epsilon \cdot b^2 \in [\mathbb{D}^*, \mathbb{D}^*] & \mbox{if } a = 0 \end{array}
                        \text{ where }\;\epsilon = \pm 1 \right\} \,.
        \end{align*}
\end{lemma}
\begin{proof}
        This is a straightforward calculation.
\end{proof}

\begin{proposition}\label{pr:Dtr}
        Let $G = \SL_2(\mathbb{D})$ and let $\tau = \mtr{0}{1}{-1}{0} \in G$.
        Let $X$ be the projective line over $\mathbb{D}$, i.e.\@ 
        $X = \{ \vct{a}{b} \mathbb{D} \neq 0 \mid a,b \in \mathbb{D} \}$.
        Then the following are equivalent: 
\begin{enumerate}
\item $C_G(\tau)$ is transitive on $X$;
        \item $a^2 + aba^{-1}b \in (\mathbb{D}^*)^2 \ [\mathbb{D}^*, \mathbb{D}^*]$ for all $a,b \in \mathbb{D}^*$;
        \item $1 + a^2 \in (\mathbb{D}^*)^2 \ [\mathbb{D}^*, \mathbb{D}^*]$ for all $a \in \mathbb{D}^*$.
\end{enumerate}
\end{proposition}
\begin{proof}
        Since $a^2 + aba^{-1}b = a^2(1 + a^{-1}ba^{-1}b)$, we have $a^2 + aba^{-1}b \in (\mathbb{D}^*)^2 \ [\mathbb{D}^*, \mathbb{D}^*]$
        if and only if $1 + a^{-1}ba^{-1}b \in (\mathbb{D}^*)^2 \ [\mathbb{D}^*, \mathbb{D}^*]$.
        Equivalence between (ii) and (iii) follows by replacing $a^{-1}b$ by $a$ in the latter term.

Assume now that (ii) holds.
        Let $a,b \in \mathbb{D}^*$ be arbitrary; we want to show that there exists some $g \in C_G(\tau)$
        mapping $\vct{a}{b}$ to $\vct{z}{0}$ for some $z \in \mathbb{D}^*$.
        By~(ii), we know that there is some $c \in \mathbb{D}^*$ such that
        $b^{-2} + b^{-1} a^{-1} b a^{-1} \equiv c^{-2} \mod{[\mathbb{D}^*, \mathbb{D}^*]}$.
        Let $g := \mtr{cb^{-1}}{ca^{-1}}{-ca^{-1}}{cb^{-1}}$.
        Then
        $\det(g) \equiv c^2 (b^{-2} + b^{-1} a^{-1} b a^{-1}) \equiv 1$,
        i.e.\@ $g \in G$.
        Moreover, $g \vct{a}{b} = \vct{z}{0}$ for $z = c(b^{-1}a + a^{-1}b)$, proving
        that $C_G(\tau)$ acts transitively on $X$.

        Conversely, assume that $C_G(\tau)$ acts transitively on $X$.
        Let $a,b \in \mathbb{D}^*$ be arbitrary; then there exists some $g \in C_G(\tau)$ mapping
        $\vct{1}{0} \mathbb{D}$ to $\vct{a}{b} \mathbb{D}$, i.e.\@ there is some $z \in \mathbb{D}^*$ such that $g$
        maps $\vct{z}{0}$ to $\vct{a}{b}$.
        By Lemma~\ref{le:DCG}, we know that $g$ has the form $g = \mtr{x}{y}{-y}{x}$
        with $x^2 + xyx^{-1}y \in [\mathbb{D}^*, \mathbb{D}^*]$.
        Then $g \vct{z}{0} = \vct{xz}{-yz}$, and hence $a = xz$ and $b = -yz$.
        Hence
        $a^2 + aba^{-1}b = xzxz + xzyx^{-1}yz = xzx^{-1} \cdot (x^2 + xyx^{-1}y) \cdot z$,
        and since $x^2 + xyx^{-1}y \in [\mathbb{D}^*, \mathbb{D}^*]$, this implies
        $a^2 + aba^{-1}b \equiv xzx^{-1}z \equiv z^2 \mod{[\mathbb{D}^*, \mathbb{D}^*]}$.
        Since $a,b \in \mathbb{D}^*$ were arbitrary, this proves~(ii).
\end{proof}

\begin{proposition}
        Let $G = \PSL_2(\mathbb{D})$, let $\tau = \mtr{0}{1}{-1}{0} \in \SL_2(\mathbb{D})$, and let $\tilde\tau$
        be the image of $\tau$ in $G$.
        Let $X$ be the projective line over $\mathbb{D}$, i.e.\@ 
        $X = \{ \vct{a}{b} \mathbb{D} \mid a,b \in \mathbb{D}, \text{not both zero} \}$.
        Then the following are equivalent: 
        \begin{enumerate}
        \item $C_G(\tilde\tau)$ is transitive on $X$;
        \item $PC_G(\tau)$ is transitive on $X$;
        \item $a^2 + aba^{-1}b \in \{ \pm 1 \} \cdot (\mathbb{D}^*)^2 \ [\mathbb{D}^*, \mathbb{D}^*]$ for all $a,b \in \mathbb{D}^*$;
        \item $1 + a^2 \in \{ \pm 1 \} \cdot (\mathbb{D}^*)^2 \ [\mathbb{D}^*, \mathbb{D}^*]$ for all $a \in \mathbb{D}^*$.
        \end{enumerate}
\end{proposition}
\begin{proof}
        The equivalence between (i) and (ii) follows immediately from the definition of the projective centralizer $PC_G(\tau)$.
        The other equivalences are shown exactly as in the proof of Proposition~\ref{pr:Dtr} above.
\end{proof}

\begin{corollary}
        \begin{itemize}
            \item[\rm (i)]
                Let $G = \SL_2(\mathbb{D})$, and assume that for all $a \in \mathbb{D}^*$, we have
                $1 + h_a \in H$.
        Then $C_G(\tau)$ acts transitively on $X$.
            \item[\rm (ii)]
                Let $G = \PSL_2(\mathbb{D})$, and assume that for all $a \in \mathbb{D}^*$, we have
                $1 + h_a \in \{ \pm 1 \} \cdot H$.
                Then $C_G(\tilde\tau)$ acts transitively on $X$.
        \end{itemize}
\end{corollary}
\begin{proof}
        We only show (i). The proof of (ii) is completely similar.
        So let $a \in \mathbb{D}^*$ be arbitrary, and assume that $1 + h_a = h \in H$.
        Then $1 + 1 h_a = 1 h$, i.e.\@ $1 + a^2 = 1 h$.
        Write $h = h_{x_1} \dotsm h_{x_n}$ with $x_1,\dots,x_n \in \mathbb{D}^*$.
        Then
        $1 h = x_n \dotsm x_1 \cdot 1 \cdot x_1 \dotsm x_n \equiv (x_1 \dotsm x_n)^2 \mod{[\mathbb{D}^*, \mathbb{D}^*]}$,
        and hence $1 + a^2 = 1h \in (\mathbb{D}^*)^2 \ [\mathbb{D}^*, \mathbb{D}^*]$.
        So (iii) of Proposition~\ref{pr:Dtr} holds, and therefore the group $C_G(\tau)$ acts transitively on $X$.
\end{proof}

A natural extension of the study of the obvious flip would be to study its close relatives, the \Defn{semi-obvious} flips, which are obtained by composing the obvious flip with a field (anti-)automorphism. 

\vspace{1cm}

\noindent Authors' addresses:

\vspace{.8cm}

\noindent
Tom De Medts \\
Department of Pure Mathematics and Computer Algebra \\
Ghent University \\
Krijgslaan 281, S22 \\
B-9000 Gent \\
Belgium \\
e-mail: {\tt tdemedts@cage.ugent.be}

\vspace{.8cm}

\noindent 
Ralf Gramlich, Max Horn \\
TU Darmstadt \\
Fachbereich Mathematik \\
Schlo\ss gartenstra\ss e 7 \\
64289 Darmstadt \\
Germany \\
e-mail: {\tt gramlich@mathematik.tu-darmstadt.de} \\
{\tt mhorn@mathematik.tu-darmstadt.de}

\vspace{.8cm}

\noindent Second author's alternative address:

\noindent University of Birmingham \\
School of Mathematics \\
Edgbaston \\
Birmingham \\
B15 2TT \\
United Kingdom \\
e-mail: {\tt ralfg@maths.bham.ac.uk}

\end{document}